\documentclass{gtart_h}   

\def\ifplaintex{\expandafter\ifx\csname documentclass\endcsname\relax}

\ifplaintex 
\else
\fi

\expandafter\ifx\csname beginpicture\endcsname\relax
\expandafter\ifx\csname documentclass\endcsname\relax
\input pictex \else
\input prepictex \input pictex \input postpictex \fi\fi

\def\gt{{\mathsurround=0pt\it $\cal G\mskip-2mu$eometry \&\ 
$\cal T\!\!$opology}}        

\def\gtp{{\mathsurround=0pt\it $\cal G\mskip-2mu$eometry \&\ 
$\cal T\!\!$opology $\cal P\!$ublications}}  

\def\lognumber#1{\def\thelognumber{#1}}
\def\volumenumber#1{\def\thevolumenumber{#1}}
\def\papernumber#1{\def\thepapernumber{#1}}
\def\volumeyear#1{\def\thevolumeyear{#1}}

\def\pagenumbers#1#2{\def\startpage{#1}\def\finishpage{#2}}
\def\published#1{\def\publishdate{#1}}
\def\proposed#1{\def\theproposer{#1}}
\def\seconded#1{\def\theseconders{#1}}
\def\received#1{\def\receiveddate{#1}}
\def\revised#1{\def\reviseddate{#1}}
\def\accepted#1{\def\accepteddate{#1}}

\long\def\asciiabstract#1{\long\def\theasciiabstract{#1}}
\def\asciikeywords#1{\def\theasciikeywords{#1}}

\let\\\par\let\thelognumber\relax
\let\thevolumenumber\relax\let\thepapernumber\relax
\let\thevolumeyear\relax\let\thesamplenumber\relax\let\startpage\relax
\let\finishpage\relax\let\publishdate\relax\let\receiveddate\relax
\let\reviseddate\relax\let\accepteddate\relax\let\theasciititle\relax
\let\theasciiauthors\relax
\let\theasciiabstract\relax\let\theasciikeywords\relax
\let\theasciiemail\relax\let\theshortauthors\relax\let\theshorttitle\relax

\long\def\maketitlep{   

\count0=\startpage

\gt\hfill      
\beginpicture
\setcoordinatesystem units <0.33truein, 0.33truein> point at 2.2 0.9
\setplotsymbol ({$\cal G$})
\plotsymbolspacing=9truept
\circulararc 315 degrees from 0 1 center at 0 0
\setplotsymbol ({$\cal T$})
\circulararc 315 degrees from 1 -1 center at 1 0
\endpicture
\break
{\small\ifx\thesamplenumber\relax 
Volume \else Sample
\fi\thevolumenumber\ (\thevolumeyear)
\startpage--\finishpage\nl
Published: \publishdate}
\vglue 0.5truein plus 0.4fil minus 0.1truein

{\parskip=0pt\leftskip 0pt plus 1fil\def\\{\par\smallskip}{\ifplaintex\large
\else\Large\fi\bf\thetitle}\par\medskip}   

\vglue 0pt plus 0.1fil 

{\parskip=0pt\leftskip 0pt plus 1fil\def\\{\par}{\sc\theauthors}
\par\medskip}

\vglue 0pt plus 0.1fil 

{\small\parskip=0pt\let\newline\\
{\leftskip 0pt plus 1fil\def\\{\par}{\sl\theaddress}\par}
\expandafter\ifx\theemail\relax    
\relax\else\vglue 5pt plus 0.02fil minus 2pt\def\\{\stdspace{\rm 
and}\stdspace} 
\cl{Email:\stdspace\tt\theemail}\fi
\ifx\theurl\relax                  
\relax\else\vglue 5pt plus 0.02fil minus 2pt\def\\{\stdspace{\rm 
and}\stdspace}
\cl{URL:\stdspace\tt\theurl}\fi\par}

\vglue 7pt plus 0.3fil minus 3pt

{\bf Abstract}
\vglue 5pt plus 0.1fil minus 2pt

\theabstract

\vglue 7pt plus 0.3fil minus 3pt

{\bf AMS Classification numbers}\quad Primary:\quad \theprimaryclass

Secondary:\quad \thesecondaryclass

\vglue 5pt plus 0.3fil minus 2pt

{\bf Keywords}\quad \thekeywords

\vglue 10pt plus 0.5fil minus 5pt

{\small  Proposed: \theproposer\hfill Received: \receiveddate\nl
Seconded: \theseconders\hfill 
\ifx\reviseddate\relax                         
Accepted: \accepteddate                        
\else
Revised: \reviseddate                          
\fi}
\eject
}       

\let\maketitlepage\maketitlep
\let\maketitle\maketitlepage

\font\phead=cmsl9 scaled 950
\font=cmsl9 scaled 1050
\font\pnum=cmbx10 scaled 913
\font\lnum=cmbx10 
\font\pfoot=cmsl9 scaled 950
\font=cmsl9 scaled 1050
\ifplaintex
\headline{\vbox to 0pt{\vskip -4.5mm\line{\small\phead\ifnum
\count0=\startpage ISSN 1364-0380 (on line)
1465-3060 (printed) \hfill {\pnum\folio}\else\ifodd\count0\def\\{ }%
\ifx\theshorttitle\relax\thetitle\else\theshorttitle\fi\hfill{\pnum\folio}
\else\def\\{ and }{\pnum\folio}\hfill\ifx\theshortauthors\relax\theauthors
\else\theshortauthors\fi\fi\fi}\vss}}
\footline{\vbox to 0pt{\vglue 0mm\line{\small\pfoot\ifnum\count0=\startpage
\copyright\ \gtp\hfill\else
\gt, Volume \thevolumenumber\ (\thevolumeyear)\hfill\fi}\vss
}}
\else
\makeatletter
\def\@oddhead{{\small\lhead\ifnum\count0=\startpage ISSN 1364-0380 (on line)
1465-3060 (printed) \hfill {\lnum\number\count0}\else\ifodd\count0
\def\\{ }\ifx\theshorttitle\relax \thetitle \else\theshorttitle\fi\hfill
{\lnum\number\count0}\else\def\\{ and }{\lnum\number\count0}
\hfill\ifx\theshortauthors\relax 
\theauthors\else\theshortauthors\fi\fi\fi}}\def\@evenhead{\@oddhead}
\def\@oddfoot{\small\lfoot\ifnum\count0=\startpage\copyright\ \gtp\hfill\else
\gt, Volume \thevolumenumber\ (\thevolumeyear)\hfill\fi}
\def\@evenfoot{\@oddfoot}
\makeatother
\fi

\newwrite\gtoutfile
\long\gdef\makeheadfile{  
{\def\\{, }\def\s{ }
\immediate\openout\gtoutfile head.xxx
\immediate\write\gtoutfile{Proxy-for: \ifx\theasciiauthors\relax
\theauthors\else\theasciiauthors\fi\s<\ifx\theasciiemail\relax\theemail\else\theasciiemail\fi>}
\immediate\write\gtoutfile{\noexpand\\}
\immediate\write\gtoutfile{Authors: \ifx\theasciiauthors\relax
\theauthors\else\theasciiauthors\fi}
{\def\\{ }\immediate\write\gtoutfile{Title: \ifx\theasciititle\relax
\thetitle\else\theasciititle\fi}}
\immediate\write\gtoutfile{Subj-class: GT or SG or MG etc}
\immediate\write\gtoutfile{MSC-class: \theprimaryclass\ifx\thesecondaryclass\relax\else, \thesecondaryclass\fi}
\immediate\write\gtoutfile{Journal-ref: Geom. Topol. \thevolumenumber
(\thevolumeyear) \startpage-\finishpage}
\immediate\write\gtoutfile{Comments: Published by Geometry and Topology at}
\immediate\write\gtoutfile{\s\s http://www.maths.warwick.ac.uk/gt/GTVol\thevolumenumber/paper\thepapernumber.abs.html}
\immediate\write\gtoutfile{\noexpand\\}
\immediate\write\gtoutfile{}
\ifx\theasciiabstract\relax
\immediate\write\gtoutfile{\theabstract}\else
\immediate\write\gtoutfile{\theasciiabstract}\fi
\immediate\write\gtoutfile{}
\immediate\write\gtoutfile{\noexpand\\}
\immediate\write\gtoutfile{}
\immediate\closeout\gtoutfile}}  

\def\maketitlepage{\maketitlep\makeheadfile}
\let\maketitle\maketitlepage

\lognumber{343}
\volumenumber{8}\papernumber{10}\volumeyear{2004}
\pagenumbers{413}{474}
\received{23 July 2003}
\revised{17 January 2004}
\published{14 February 2004}
\accepted{14 February 2004}
\proposed{Robion Kirby}
\seconded{Shigeyuki Morita, Benson Farb}

\usepackage{amssymb,amsmath,graphicx,amscd,psfrag}             
\let\relabel\psfrag
\def\adjustrelabel <#1,#2> #3#4{%
\psfrag{#3}{\smash{\rlap{\kern #1 \raise #2\hbox{#4}}}}}

\def\Mapright#1#2{\,\smash{\mathop{\smash{\hbox to #1pt{\rightarrowfill}}
\vphantom\rightarrow}\limits^{#2}\,}}
\def\col{\hbox{$\,:\,$}} 
\def\larrow{\mathop{\longrightarrow}\limits}
\def\edge{$1$--simplex}
\def\edges{$1$--simplices}
\newcommand{\Jscr}{\mathcal J}
\newcommand{\Hscr}{\mathcal H}
\newcommand{\vol}{\operatorname{vol}}
\newcommand{\Ker}{\operatorname{Ker}}
\newcommand{\Hom}{\operatorname{Hom}}
\renewcommand{\Im}{\operatorname{Im}}

\newcommand{\cs}{\operatorname{cs}}
\newcommand{\Isom}{\operatorname{Isom}}
\newcommand{\PSL}{\operatorname{PSL}}
\newcommand{\PGL}{\operatorname{PGL}}
\newcommand{\SL}{\operatorname{SL}}
\newcommand\calE{{\mathcal E}}
\newcommand\calR{{\mathcal R}}
\newcommand\eeprebloch{{\calE\prebloch}}
\newcommand\eebloch{{\calE\bloch}}
\newcommand\prebloch{{\mathcal P}}
\newcommand\bloch{{\mathcal B}}
\newcommand\FT{\operatorname{FT}}
\newcommand\liftFT{\widehat\FT}
\newcommand\eprebloch{{\widehat\prebloch}}
\newcommand\ebloch{{\widehat\bloch}}
\newcommand\C{\mathbb C}
\newcommand\R{\mathbb R}
\newcommand\Q{\mathbb Q}
\newcommand\CP{\C P}
\newcommand\Z{\mathbb Z}
\newcommand\bfw{{\bf w}}
\renewcommand\H{\mathbb H}
\newcommand\discrete{^\delta}

\newcommand\themap{\lambda}
\newtheorem{theorem}{Theorem}[section]
\newtheorem{lemma}[theorem]{Lemma}
\newtheorem{proposition}[theorem]{Proposition}
\newtheorem{corollary}[theorem]{Corollary}
\newtheorem*{scholium*}{Scholium}
\theoremstyle{definition}
\newtheorem{definition}[theorem]{Definition}
\newtheorem{remark}[theorem]{Remark}
\newtheorem*{remark*}{Remark}
\newtheorem*{example*}{Example}
\newtheorem*{assumption*}{Assumption}
\newtheorem{example}[theorem]{Example}
\newcommand\Cover{\widehat\C}
\begin{document}
\title{Extended Bloch group and the\\Cheeger--Chern--Simons class}
\author{Walter D Neumann} 
\address{Department of Mathematics, Barnard College\\Columbia
  University, New York, NY 10027, USA} 
\email{neumann@math.columbia.edu}

\begin{abstract}
  We define an extended Bloch group and show it is naturally
  isomorphic to $H_3(\PSL(2,\C)^\delta;\Z)$. Using the Rogers
  dilogarithm function this leads to an exact simplicial formula for
  the universal Cheeger--Chern--Simons class on this homology group.  It also
  leads to an independent proof of the analytic relationship between
  volume and Chern--Simons invariant of hyperbolic $3$--manifolds
  conjectured in \cite{neumann-zagier} and proved in \cite{yoshida},
  as well as effective formulae for the Chern--Simons invariant of a
  hyperbolic $3$--manifold.
\end{abstract}
\asciiabstract{We define an extended Bloch group and show it is
  naturally isomorphic to H_3(PSL(2,C)^\delta;Z).  Using the Rogers
  dilogarithm function this leads to an exact simplicial formula for
  the universal Cheeger-Chern-Simons class on this homology group.  It
  also leads to an independent proof of the analytic relationship
  between volume and Chern-Simons invariant of hyperbolic 3-manifolds
  conjectured by Neumann and Zagier [Topology 1985] and proved by
  Yoshida [Invent. Math. 1985] as well as effective formulae for
  the Chern-Simons invariant of a hyperbolic 3-manifold.}

\asciikeywords{Extended Bloch group, Cheeger-Chern-Simons class, hyperbolic,
3-manifold}

\maketitlepage

\section{Introduction} 
There are several variations of the definition of the Bloch group in
the literature; by \cite{dupont-sah} they differ at most by torsion
and they agree with each other for algebraically closed fields.  In
this paper we shall use the following.

\begin{definition}\label{def-bloch}
Let $k$ be a field.  The \emph{pre-Bloch group $\prebloch(k)$} is the
quotient of the free $\Z$--module $\Z (k-\{0,1\})$ by all instances of
the following relation:
\begin{equation}
[x]-[y]+[\frac yx]-[\frac{1-x^{-1}}{1-y^{-1}}]+[\frac{1-x}{1-y}]=0,
\label{fiveterm}
\end{equation}
This relation is usually called the \emph{five term relation}.  The
\emph{Bloch group $\bloch(k)$} is the kernel of the map
$$
\prebloch(k)\to k^*\wedge_\Z k^*,\quad [z]\mapsto 2\bigl(z\wedge(1-z)\bigr).
$$
\end{definition}

(In \cite{neumann-yang} the additional relations
$$
[x]=[1-\frac 1x]=[\frac 1{1-x}]=-[\frac1x]=-[\frac{x-1}x]=-[1-x]
\label{invsim}
$$
were used.  These follow from the five term relation when $k$ is
algebraically closed, as shown by
Dupont and Sah \cite{dupont-sah}.  Dupont and Sah use a different five
term relation but it is conjugate to the one used here by
$z\mapsto\frac1z$.)

There is an exact sequence due to Bloch and Wigner: 
$$
0\to \Q/\Z\to H_3(\PSL(2,\C)\discrete;\Z)\to \bloch(\C)\to 0.
$$
The superscript $\delta$ means ``with discrete topology.''
\renewcommand\discrete{} We will omit it from now on.

$\bloch(\C)$ is known to be uniquely divisible, so it has canonically
the structure of a $\Q$--vector space (Suslin \cite{suslin}). It's
$\Q$--dimension is infinite and conjectured to be countable (the
``Rigidity Conjecture,'' equivalent to the conjecture that
$\bloch(\C)=\bloch(\overline\Q)$, where $\overline\Q$ is the field of
algebraic numbers).  In particular, the $\Q/\Z$ in the Bloch--Wigner
exact sequence is precisely the torsion of $H_3(\PSL(2,\C);\Z)$, so
any finite torsion subgroup is cyclic.

\medskip In the present paper we define an \emph{extended Bloch group}
$\ebloch(\C)$ by replacing $\C-\{0,1\}$ in the definition of
$\bloch(\C)$ by an abelian cover and appropriately lifting the five
term relation (\ref{fiveterm}). Our main results are that we can lift
the Bloch--Wigner map $H_3(\PSL(2,\C)\discrete;\Z)\to \bloch(\C)$ to an
isomorphism
$$
\themap\co H_3(\PSL(2,\C)\discrete;\Z)\stackrel\cong\longrightarrow
\ebloch(\C)
$$
Moreover, the ``Roger's dilogarithm function'' (see below) gives a
natural map
$$R\co\ebloch(\C)\to\C/\pi^2\Z.$$ We
show that the composition $$R\circ\themap\co
  H_3(\PSL(2,\C)\discrete;\Z)\to\C/\pi^2\Z$$ is the
Cheeger--Chern--Simons class (cf \cite{cheeger-simons}), so it can also
be described as $i(\vol + i \cs)$, where $cs$ is the universal
Chern--Simons class.  It has been a longstanding problem to provide
such a computation of the Chern--Simons class. Dupont in
\cite{dupont1} gave an answer modulo $\pi^2\Q$ and our computation is
a natural lift of his.

We show that any complete hyperbolic 3--manifold $M$ of finite volume
has a natural ``fundamental class'' in $H_3(\PSL(2,\C)\discrete;\Z)$.
For compact $M$ the existence of this class is easy and well known:
$M=\H^3/\Gamma$ is a $K(\Gamma,1)$--space, so the inclusion $\Gamma\to
\Isom^+(\H^3) = \PSL(2,\C)$ induces $H_3(M)=H_3(\Gamma)\to
H_3(\PSL(2,\C))$, and the class in question is the image of the
fundamental class $[M]\in H_3(M)$.  For non-compact $M$ 
the existence is shown in Section
\ref{sec:3-manifolds}. The Cheeger--Chern--Simons class applied to this
class gives $i(\vol(M)+i\cs(M))$ where $\vol(M)$ is hyperbolic volume
and $\cs(M)$ is the Chern--Simons invariant. This Chern--Simons
invariant is defined modulo $\pi^2\Z$ and for compact $M$ it arises
here as the Chern--Simons invariant for the flat $PSL(2,\C)$--bundle
over $M$.  According to J.~Dupont, this is the same as the
Chern--Simons invariant for the Riemannian connection over $M$ (private
communication; this is proved modulo $6$--torsion in
\cite{dupont-kamber} and our results also provide confirmation---see
\cite{coulson-et-al} for discussion).
However, the Chern--Simons invariant for the Riemannian connection
carries slightly more information, since it is defined modulo $2\pi^2$
rather than modulo $\pi^2$. For non-compact $M$ we show that the
Chern--Simons invariant is the one defined by Meyerhoff
\cite{meyerhoff}, which is naturally only defined modulo $\pi^2$.

This fundamental class of $M$ in $H_3(\PSL(2,\C);\Z)$ determines an
element $\hat\beta(M)$\break $\in\ebloch(\C)$.  Our results describe
$\hat\beta(M)$ directly in terms of an ideal triangulation of $M$, as
a lift of the Bloch invariant $\beta(M)\in\bloch(\C)$ defined in
\cite{neumann-yang}.  We need a ``true'' ideal triangulation, which is
more restrictive than the ``degree 1'' ideal triangulations that
sufficed in \cite{neumann-yang}.  The ideal triangulations resulting
from Dehn filling that are used by the programs Snappea \cite{snappea}
and Snap \cite{snap} (see also \cite{coulson-et-al}) are not true, but
we describe $\hat\beta(M)$ in terms of these ``Dehn filling
triangulations'' in Theorems \ref{th:conj true} and \ref{th:deg one},
and deduce practical simplicial formulae for Chern--Simons invariants
of hyperbolic 3--manifolds.

The formula of Theorem \ref{th:conj true} exhibits directly the
analytic relationship, conjectured in \cite{neumann-zagier} and proved
in \cite{yoshida}, between volume and Chern--Simons invariant of
hyperbolic manifolds, and gives an independent proof of it.  A similar
formula for Chern--Simons invariant was derived from
\cite{neumann-zagier,yoshida} in \cite{neumann}, but that formula was
only accurate up to an unknown constant on any given Dehn surgery
space (the constant was conjectured to be a multiple of $\pi^2/6$, and
our results now confirm this). Snappea and versions of Snap prior to
1.10.2 use versions of that formula, and use a bootstrapping procedure
to discover the constant and obtain the Chern--Simons invariant for
many manifolds.  Thanks to Oliver Goodman, Snap now implements the
formula of Theorem \ref{th:deg one} and can thus compute the
Chern--Simons invariant for any hyperbolic manifold.

More generally, any flat $\PSL(2,\C)$ bundle over a closed oriented
$3$--manifold $M$ determines a class in $H_3(\PSL(2,\C);\Z)$ and our
results give a simplicial computation of this class as an element of
$\ebloch(\C)$. 

The main results of this paper were announced in
\cite{neumann-hilbert} and partial proofs were in the preliminary
preprint \cite{extblochorig}.
This paper corrects the tentative statement in \cite{neumann-hilbert}
(also in \cite{extblochorig}) that our map $\lambda$ is not an
isomorphism but has a cyclic kernel of order $2$.  There is a related
error in Section 6B of \cite{coulson-et-al}: it is wrongly stated that
the Cheeger--Chern--Simons class takes values in $\C/2\pi^2\Z$ rather
than $\C/\pi^2\Z$ (and therefore that the fundamental class in
$H_3(\PSL(2,\C);\Z)$ of a cusped hyperbolic $3$--manifold has a mod 2
ambiguity).

\medskip\noindent{\bf Acknowledgments}\qua The definition of the extended
Bloch group was suggested by an idea of Jun Yang, to whom I am
grateful also for many useful conversations.  I'm grateful to
Benedetti and Baseilhac, whose use in \cite{benedetti} of the
preliminary preprint \cite{extblochorig} and subsequent correspondence
led me to finish writing this paper.  This research was supported by
the Australian Research Council and by the NSF.

\section{The extended Bloch group}

We shall need a $\Z\times\Z$ cover $\Cover$ of $\C-\{0,1\}$ which can
be constructed as follows.  Let $P$ be $\C-\{0,1\}$ split along the
rays $(-\infty,0)$ and $(1,\infty)$.  Thus each real number $r$
outside the interval $[0,1]$ occurs twice in $P$, once in the upper
half plane of $\C$ and once in the lower half plane of $\C$.  We
denote these two occurrences of $r$ by $r+0i$ and $r-0i$.  We construct
$\Cover $ as an identification space from $P\times\Z\times\Z$ by
identifying
$$\begin{aligned}
(x+0i, p,q)&\sim (x-0i,p+2,q)\quad\hbox{for each }x\in(-\infty,0)\\
(x+0i, p,q)&\sim (x-0i,p,q+2)\quad\hbox{for each }x\in(1,\infty).
\end{aligned}$$
We will denote the equivalence class of $(z,p,q)$ by $(z;p,q)$.
$\Cover $ has four components:
$$
\Cover =X_{00}\cup X_{01}\cup X_{10}\cup X_{11}
$$
where $X_{\epsilon_0\epsilon_1}$ is the set of $(z;p,q)\in \Cover $ with
$p\equiv\epsilon_0$ and $q\equiv\epsilon_1$ (mod $2$).

We may think of $X_{00}$ as the Riemann surface for the multivalued function
$$\C-\{0,1\}\to\C^2,\qquad z\mapsto \bigl(\log z, -\log (1-z)\bigr).$$
Taking the branch $(\log z + 2p\pi i, -\log (1-z) + 2q\pi i)$ of this
function on the portion $P\times\{(2p,2q)\}$ of $X_{00}$ for each
$p,q\in\Z$ defines an analytic function from $X_{00}$ to $\C^2$.  In
the same way, we may think of $\Cover $ as the Riemann surface for the
multivalued function $(\log z + p\pi i, -\log
(1-z) + q\pi i)$, $p,q\in\Z$, on $\C-\{0,1\}$.

\medbreak
Consider the set
$$
\FT:=\biggl\{\biggl(x,y,\frac yx,\frac{1-x^{-1}}{1-y^{-1}},
\frac{1-x}{1-y}\biggr):x\ne y,
x,y\in\C-\{0,1\}\biggr\}\subset(\C-\{0,1\})^5
$$
of 5--tuples involved in the five term relation (\ref{fiveterm}).  
An elementary computation shows:

\begin{lemma}\label{lem:FT}
The subset $\FT^+$ of $(x_0,\dots,x_4)\in \FT$ with each $x_i$ 
in the upper half plane of $\C$ is the set of elements of $\FT$ for
which $y$ is in the upper half plane of $\C$ and $x$ is inside the
triangle with vertices $0,1,y$. Thus $\FT^+$ is connected (even
contractible).  \qed\end{lemma}

\begin{definition}\label{liftFT}
Let $V\subset(\Z\times\Z)^5$ be the subspace
\begin{multline*}
V:=\{\bigl((p_0,q_0),(p_1,q_1),(p_1-p_0,q_2),
(p_1-p_0+q_1-q_0,q_2-q_1),\\ (q_1-q_0,q_2-q_1-p_0)\bigr):p_0,p_1,q_0,q_1,q_2\in\Z\}.
\end{multline*}
Let $\liftFT_0$ denote the unique component of the inverse image of $\FT$ in
$\Cover^5$ which includes the points $\bigl((x_0;0,0),\ldots,(x_4;0,0)\bigr)$
with $(x_0,\ldots,x_4)\in \FT^+$, and define
$$
\liftFT:=\liftFT_0+V=\{
{\mathbf x}+{\mathbf v}:{\mathbf x}\in
\liftFT_0\text{ and }{\mathbf v}\in V\},
$$
where we are using addition to denote the action of $(\Z\times\Z)^5$
by covering transformations on $\Cover^5$. (Although we do not need
it, one can show that the action of $2V$ takes $\liftFT_0$ to itself,
so $\liftFT$ has $2^5$ components, determined by the parities of
$p_0,p_1,q_0,q_1,q_2$.)

Define $\eprebloch(\C)$ as the free
$\Z$--module on $\Cover$ factored by all instances of the relations:%
\begin{equation}\label{five term}
\sum_{i=0}^4(-1)^i(x_i;p_i,q_i)=0\quad
\text
{with $\bigl((x_0;p_0,q_0),\dots,(x_4;p_4,q_4)\bigr)\in\liftFT$}
\end{equation}
and
\begin{equation}
(x;p,q)+(x;p',q')=(x;p,q')+(x;p',q)\label{transfer}
\quad\text{with $p,q,p',q'\in\Z$}
\end{equation}
We shall denote the class of $(z;p,q)$ in
$\eprebloch(\C)$ by $[z;p,q]$
\end{definition}

We call relation (\ref{five term}) the \emph{lifted five term relation}.
We shall see that its precise form arises naturally in several
contexts. In particular, we give it a geometric interpretation in
Section~\ref{simplex parameters}.  

We call relation (\ref{transfer}) the \emph{transfer relation}.  It is
almost a consequence of the lifted five term relation, since we shall
see that the effect of omitting it would be to replace
$\eprebloch(\C)$ by $\eprebloch(\C)\oplus\Z/2$, with $\Z/2$ generated
by an element $\kappa:=[x;1,1]+[x;0,0]-[x;1,0]-[x;0,1]$ which is
independent of $x$.

\begin{lemma}\label{le:dehn}
There is a well-defined homomorphism
$$\nu\co\eprebloch(\C)\to \C\wedge_\Z\C
$$ defined on generators by
$[z;p,q]\mapsto (\log z + p\pi i)\wedge (-\log(1-z) +q\pi i)$.  
\end{lemma}
\begin{proof}
We must verify that $\nu$ vanishes on the relations that define
$\eprebloch(\C)$. This is trivial for the transfer relation
(\ref{transfer}).
We shall show that the lifted five term
relation is the most general lift of the five term relation
(\ref{fiveterm}) for which $\nu$ vanishes.  If one applies $\nu$ to an
element $\sum_{i=0}^4(-1)^i[x_i;p_i,q_i]$ with
$(x_0,\ldots,x_4)=(x,y,\ldots)\in \FT^+$ one obtains after simplification:
\begin{multline*}
\bigl((q_0-p_2-q_2+p_3+q_3)\log x
+(p_0-q_3+q_4)\log(1-x)+(-q_1+q_2-q_3)\log y +{}\\ {}+
(-p_1+p_3+q_3-p_4-q_4)\log(1-y)+(p_2-p_3+p_4)\log(x-y)\bigr)\wedge \pi i.
\end{multline*}
An elementary linear algebra computation shows that this vanishes
identically if and only if $p_2=p_1-p_0$, $p_3=p_1-p_0+q_1-q_0$,
$q_3=q_2-q_1$, $p_4=q_1-q_0$, and $q_4 = q_2-q_1-p_0$, as in the
lifted five term relation.  The vanishing of $\nu$ for the general
lifted five term relation now follows by analytic continuation.
\end{proof}
\begin{definition}
Define $\ebloch(\C)$ as the kernel of
$\nu\co\eprebloch(\C)\to\C\wedge\C$.
\end{definition}
 
Define
$$
R(z;p,q)={\calR}(z)+\frac{\pi i}2(p\log(1-z)+q\log z) -\frac{\pi^2}6
$$ 
where $\calR$ is the Rogers dilogarithm function
$$
{\calR}(z)=\frac12\log(z)\log(1-z)-\int_0^z\frac{\log(1-t)}tdt.$$
Then we have:

\begin{proposition}\label{R}
$R$ gives a well defined map $R\co \Cover\to \C/\pi^2\Z$.  The relations
which define $\eprebloch(\C)$ are functional equations for $R$ modulo
$\pi^2$ (the lifted five term relation is in fact the most general
lift of the five term relation (\ref{fiveterm}) with this property).
Thus $R$ also gives a homomorphism
$R\co\eprebloch(\C)\to\C/\pi^2\Z$.
\end{proposition}

\begin{proof}
  If one follows a closed path from $z$ that goes anti-clockwise
  around the origin it is easily verified that $R(z;p,q)$ is replaced
  by $R(z;p,q)+\pi i \log(1-z) -q\pi^2=R(z;p+2,q)-q\pi^2$.  Similarly,
  following a closed path clockwise around $1$ replaces $R(z;p,q)$ by
  $R(z;p,q+2)+p\pi^2$.  Thus $R$ modulo $\pi^2$ is well defined on
  $\Cover$ (and $R$ modulo $2\pi^2$ is well defined on $X_{00}$; in
  fact $R$ itself is well defined on a $\Z$ cover of $X_{00}$ which is
  the universal nilpotent cover of $\C-\{0,1\}$).

It is well known that
$L(z):=\calR(z)-\frac{\pi^2}6$ satisfies the functional equation
$$
L(x)-L(y)+L\bigl(\frac yx\bigr)-L\bigl(\frac{1-x^{-1}}{1-y^{-1}}\bigr)+
L\bigl(\frac{1-x}{1-y}\bigr)=0$$ 
for
$0<y<x<1$.  Since the 5--tuples involved in this equation are on the
boundary of $\FT^+$, the functional equation
$$\sum(-1)^iR(x_i;0,0)=0$$
is valid by analytic continuation on the whole of $\FT^+$.  Now
$$\sum(-1)^iR(x_i;p_i,q_i)$$
differs from this by 
$$\frac{\pi i}2\sum(-1)^i(p_i\log(1-x_i)+q_i\log x_i).
$$ 
For $(x_0,\dots,x_4)=(x,y,\dots)\in \FT^+$ this equals
\begin{multline*}
\frac{\pi i}{2}\bigl((q_0-p_2-q_2+p_3+q_3)\log x
+(p_0-q_3+q_4)\log(1-x)+(-q_1+q_2-q_3)\log y +{}\\ {}+
(-p_1+p_3+q_3-p_4-q_4)\log(1-y)+(p_2-p_3+p_4)\log(x-y)\bigr)
\end{multline*}
and as in the proof of Lemma \ref{le:dehn}, this vanishes
identically on $\FT^+$ if and only if the $p_i$
and $q_i$ are as in the lifted five term relation.  Thus the lifted
five-term relation gives a functional equation for $R$ when
$(x_0,\ldots,x_4)\in \FT^+$.  By analytic continuation, it is a
functional equation for $R$ mod $\pi^2$ in general.  The transfer
relation is trivially a functional equation for $R$.
\end{proof}

Our main result is the following:
\begin{theorem}\label{th:main}
  There exists an isomorphism $\themap\co
  H_3(\PSL(2,\C)\discrete;\Z)\to\ebloch(\C)$ such that the composition
  $\themap\circ R\co H_3(\PSL(2,\C)\discrete;\Z)\to\C/\pi^2\Z$ is
  the characteristic class given by $i(\vol + i \cs)$.
\end{theorem}

To describe the map $\themap$ we must
give a geometric interpretation of $\Cover$.
 
\section{Parameters for ideal hyperbolic simplices}
\label{simplex parameters}

In this section we shall interpret $\Cover$ as a space of parameters for
what we call ``combinatorial flattenings'' of ideal hyperbolic
simplices.  We need this to define the above map $\themap$.  It also
gives a geometric interpretation of the lifted five term relation.

We shall denote the standard compactification of $\H^3$ by $\overline
\H^3 = \H^3\cup\CP^1$. An ideal simplex $\Delta$ with vertices
$z_0,z_1,z_2,z_3\in\CP^1$ is determined up to congruence by the
cross-ratio
$$
z=[z_0\col z_1\col z_2\col z_3]=\frac{(z_2-z_1)(z_3-z_0)}{(z_2-z_0)(z_3-z_1)}.
$$
Permuting the vertices by an even (ie,
orientation preserving) permutation replaces $z$ by one of
$$
z,\quad z'=\frac 1{1-z}, \quad\text{or}\quad z''=1-\frac 1z.
$$
The parameter $z$ lies in the upper half plane of $\C$ if the
orientation induced by the given ordering of the vertices agrees with
the orientation of $\H^3$.  But we allow simplices whose vertex
ordering does not agree with their orientation. We also allow
degenerate ideal simplices whose vertices lie in one plane, so the
parameter $z$ is real.  However, we always require that the vertices
are distinct.  Thus the parameter $z$ of the simplex lies in
$\C-\{0,1\}$ and every such $z$ corresponds to an ideal simplex.

\begin{figure}[ht!]\small\nocolon
\relabel{1}{$ $}
\relabel{2}{$ $}
\relabel{3}{$ $}
\relabel{4}{$ $}
\relabel{z}{$z$}
\relabel{zz}{$z$}
\relabel{z1}{$z_0$}
\relabel{z2}{$z_1$}
\relabel{z3}{$z_2$}
\adjustrelabel <1pt,0pt> {z4}{$z_3$}
\adjustrelabel <0pt,-2pt> {zd}{$z'$}
\adjustrelabel <0pt,-1.5pt> {zzd}{$z'$}
\relabel{zdd}{$z''$}
\relabel{zzdd}{$z''$}
\centerline{
\includegraphics[width=.35\hsize]{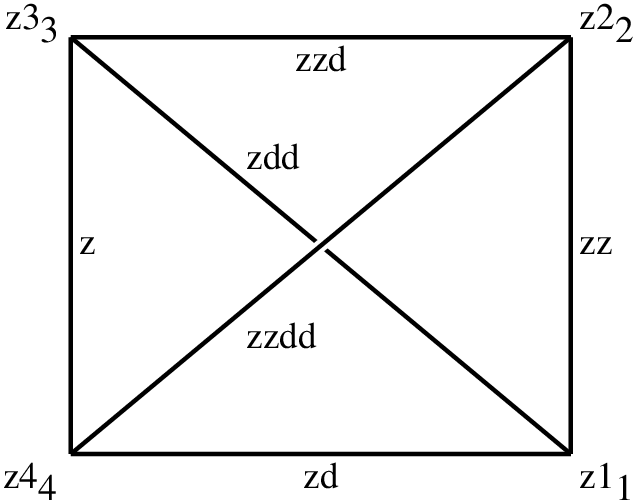}}
\caption{}\label{fig:1}
\end{figure}

There is another way of describing the cross-ratio parameter $z=
[z_0\col z_1\col z_2\col z_3]$ of a simplex. The group of orientation
preserving isometries of $\H^3$ fixing the points $z_0$ and $z_1$ is
isomorphic to $\C^*$ and the element of this $\C^*$ that takes $z_2$
to $z_3$ is $z$ (equivalently: if we position $z_0,z_1$ at $0,\infty$
in the upper half-space model of $\H^3$ then $z_3=zz_2$). Thus the
cross-ratio parameter $z$ is associated with the edge $z_0z_1$ of the
simplex\footnote{It is associated with the \emph{unoriented} edge
  since there is an orientation preserving symmetry of the simplex
  exchanging $z_0$ and $z_1$ (it also exchanges $z_2$ and $z_3$ and acts
  on $\C^*$ by $z\mapsto 1/z$).}.  The parameter associated in this
way with the other two edges $z_0z_3$ and $z_0z_2$ out of $z_0$ are
$z'$ and $z''$ respectively, while the edges $z_2z_3$, $z_1z_2$, and
$z_1z_3$ have the same parameters $z$, $z'$, and $z''$ as their
opposite edges.  See Figure~\ref{fig:1}.

Note that $zz'z''=-1$, so the sum
$$\log z + \log z' + \log z''$$ 
is an odd multiple of $\pi i$, depending on the branches of $\log$
used.  In fact, if we use standard branch of log then this sum is $\pi
i$ or $-\pi i$ depending on whether $z$ is in the upper or lower half
plane. 

The imaginary parts of $\log z$, $\log z'$, and $\log z''$ are the
dihedral angles of the ideal simplex (resp.\ their negatives if
the vertex ordering does not agree with the orientation the
simplex inherits from $\H^3$). We now consider certain ``adjustments''
of these angles by multiples of $\pi$.
\begin{definition}\label{flattening}
We shall call any triple of the form
$$
\bfw=(w_0,w_1,w_2)=(\log z +p\pi i, \log z' + q\pi i, \log z'' + r\pi i)
$$ 
with
$$p,q,r\in \Z\quad\text{and}\quad
w_0+w_1+w_2=0
$$
a \emph{combinatorial flattening} for our simplex. 

Each edge $E$ of $\Delta$ is assigned one of the
components $w_i$ of $\bfw$, with opposite edges being assigned the
same component. We call
$w_i$ the \emph{log-parameter} for the edge $E$ and denote it
$l_E(\Delta,\bfw)$. 
\end{definition}
This combinatorial flattening can be written 
$$\ell(z;p,q):=(\log z + p\pi i, -\log(1-z) + q\pi i,
\log(1-z)-\log z - (p+q)\pi i),
$$ 
and $\ell$ is then a map
of $\Cover$ to the set of combinatorial flattenings of simplices.

\begin{lemma}\label{Cover is flattenings} 
This map $\ell$ is a bijection, so $\Cover$ may be identified with
the set of all combinatorial flattenings of ideal tetrahedra.
\end{lemma}

\begin{proof}
We must show that we can recover $(z;p,q)$ from $(w_0,w_1,w_2){=}
\ell(z;p,q)$.  It clearly suffices to recover $z$. But $z=\pm e^{w_0}$
and $1-z=\pm e^{-w_1}$, and the knowledge of both $z$ and $1-z$ up to
sign determines $z$.
\end{proof}

We can give a geometric interpretation of the choice of parameters in
the five term relation (\ref{five term}).  If $z_0,\ldots,z_4$ are
five distinct points of $\partial\overline\H^3$, then each choice of
four of five points $z_0,\dots,z_4$ gives an ideal simplex. We denote
the simplex which omits vertex $z_i$ by $\Delta_i$. The cross ratio
parameters $x_i=[z_0\col \ldots\col \hat{z_i}\col \ldots\col z_4]$ of
these simplices can be expressed in terms of $x:=x_0$ and $y:=x_1$ as
follows:
$$\begin{aligned}
x_0=[z_1\col z_2\col z_3\col z_4]&=:x\\
x_1=[z_0\col z_2\col z_3\col z_4]&=:y\\
x_2=[z_0\col z_1\col z_3\col z_4]&=\frac yx\\
x_3=[z_0\col z_1\col z_2\col z_4]&=\frac {1-x^{-1}}{1-y^{-1}}\\
x_4=[z_0\col z_1\col z_2\col z_3]&=\frac {1-x}{1-y}
\end{aligned}$$
The lifted
five term relation has the form \def\term#1{(x_#1;p_#1,q_#1)}
\begin{equation}\label{general 5-term}
\sum_{i=0}^4(-1)^i \term i=0
\end{equation}
with certain relations on the $p_i$ and $q_i$.  We will give a
geometric interpretation of
these relations.

Using the map of Lemma
\ref{Cover is flattenings}, each summand in this relation
(\ref{general 5-term}) represents a choice $\ell\term i$ of combinatorial
flattening for one of the five ideal simplices. For each \edge{} $E$
connecting two of the points $z_i$ we get a corresponding linear
combination
\begin{equation}\label{edge sums}
\sum_{i=0}^4(-1)^il_E(\Delta_i,\ell\term i)
\end{equation}
of log-parameters (Definition \ref{flattening}), where we put
$l_E(\Delta_i,\ell\term i)=0$ if the line $E$ is not an edge of
$\Delta_i$.  This linear combination has just three non-zero terms
corresponding to the three simplices that meet at the edge $E$. One
easily checks that the real part is zero and the imaginary part can be
interpreted (with care about orientations) as the sum of the
``adjusted angles'' of the three flattened simplices meeting at $E$.

\begin{definition}
We say that the $\term i$ satisfy the \emph{flattening
condition} if each of the above linear combinations (\ref{edge sums})
of log-parameters is equal to zero. That is, the adjusted angle sum of
the three simplices meeting at each edge is zero.
\end{definition}

\begin{lemma}\label{geometric five term}
Relation (\ref{general 5-term}) is an instance of the lifted five term
relation (\ref{five term}) if and only if the $\term i$ satisfy
the flattening condition.
\end{lemma}
\begin{proof}
We first consider the case that $(x_0,\ldots,x_4)\in \FT^+$.  Recall
this means that each $x_i$ is in $\H$.  Geometrically, this implies
that each of the above five tetrahedra is positively oriented by the
ordering of its vertices.  
\begin{figure}[b!]\small\nocolon
\relabel{0}{$ $}
\relabel{1}{$ $}
\relabel{2}{$ $}
\relabel{3}{$ $}
\relabel{4}{$ $}
\relabel{z0}{$z_0$}
\adjustrelabel <9pt, 0pt> {z1}{$z_1$}
\adjustrelabel <3pt, 0pt> {z2}{$z_2$}
\adjustrelabel <5pt, 0pt> {z3}{$z_3$}
\adjustrelabel <0pt, 5pt> {z4}{$z_4$}
\centerline{
\includegraphics[width=.4\hsize]{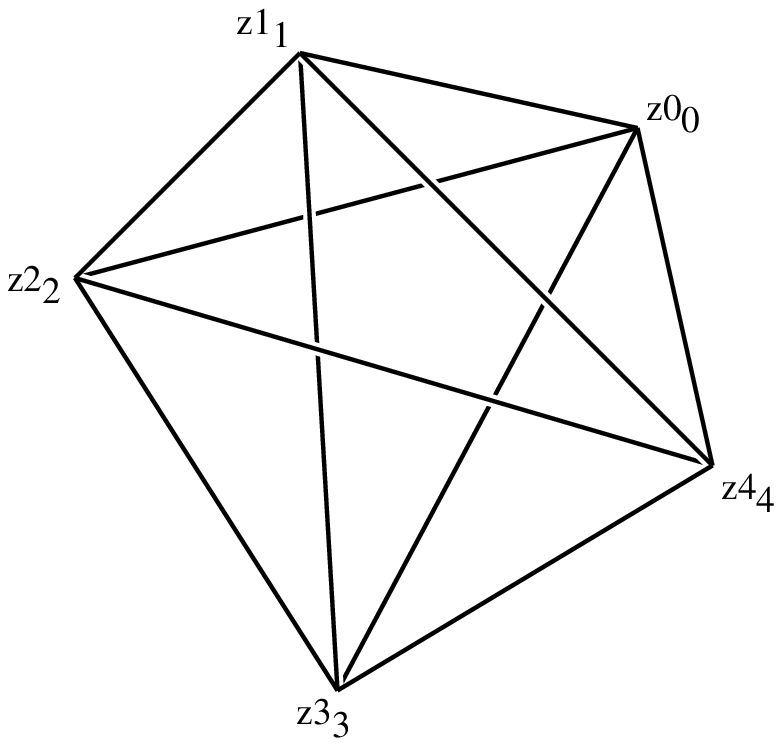}}\vspace{-6mm}
\caption{}\label{fig:2}
\end{figure}
This implies the configuration of Figure~\ref{fig:2}
with $z_1$ and $z_3$ on opposite sides of the plane of the triangle
$z_0z_2z_4$ and the line from $z_1$ to $z_3$ passing through the
interior of this triangle.  Denote the combinatorial flattening of the
$i^{th}$ simplex by $\ell(x_i;p_i,q_i)$.  If we consider the
log-parameters at the edge $z_3z_4$ for example, they are $\log x +
p_0\pi i$, $\log y + p_1 \pi i$, and $\log(y/x) + p_2\pi i$ and the
condition is that $(\log x + p_0\pi i)-(\log y + p_1
\pi i) +(\log(y/x) +p_2\pi i)=0$.  
This implies $p_2=p_1-p_0$.  Similarly the other edges lead to 
other relations among the $p_i$ and $q_i$, namely:

\begin{tabular}
{l r @{\,\,} c @{\,\,} l @{\qquad\qquad} l  r @{\,\,} c @{\,\,} l}
$z_0z_1$:&$ p_2-p_3+p_4$&$=$&$0$&$z_0z_2$:&$ -p_1+p_3+q_3-p_4-q_4$&$=$&$0$\\
$z_1z_2$:&$ p_0-q_3+q_4$&$=$&$0$&$z_1z_3$:&$ -p_0-q_0+q_2-p_4-q_4$&$=$&$0$\\
$z_2z_3$:&$ q_0-q_1+p_4$&$=$&$0$&$z_2z_4$:&$ -p_0-q_0+p_1+q_1-p_3$&$=$&$0$\\
$z_3z_4$:&$ p_0-p_1+p_2$&$=$&$0$&$z_3z_0$:&$ p_1+q_1-p_2-q_2+q_4$&$=$&$0$\\
$z_4z_0$:&$ -q_1+q_2-q_3$&$=$&$0$&$z_4z_1$:&$ q_0-p_2-q_2+p_3+q_3$&$=$&$0$.
\end{tabular}

\medskip
Elementary linear algebra verifies that these relations are equivalent
to the equations $p_2=p_1-p_0$, $p_3=p_1-p_0+q_1-q_0$, $q_3=q_2-q_1$,
$p_4=q_1-q_0$, and $q_4 = q_2-q_1-p_0$, as in the lifted five term
relation (\ref{five term}).  The lemma thus follows for
$(x_0,\ldots,x_4)\in \FT^+$. It is then true in general by analytic
continuation.
\end{proof} 
We mention here a lemma that will be useful later.
\begin{lemma}\label{lem:extend}
  If the flattenings $(x_i,p_i,q_i)$ are specified for a subset of the
  above five ideal simplices and the sum of adjusted angles is zero
  around each edge $E$ that lies on three simplices of this subset,
  then one can specify flattenings on the remaining simplices so that
  the flattening condition holds (sum of adjusted angles is zero for all
  edges $E$).  Moreover, once the flattenings are specified on three of
  the simplices, the flattenings on the final two are uniquely
  determined.
\end{lemma}
\begin{proof}
  The lemma does not depend on the ordering of $z_0, \dots, z_4$ since
  the flattening condition at an edge is purely geometric. We assume
  therefore that the specified flattenings are $(x_i;p_i,q_i)$ for
  $i=0,\dots,k$ with $0\le k < 4$.  The above proof shows that if
  $k=0$ or $1$ there is no restriction on the flattenings, while if
  $k=2$ the the condition at the edge $z_3z_4$ implies
  $p_2=p_1-p_0$ if $(x_0,\dots,x_4)\in\FT^+$, or the appropriate
  analytic continuation of this in general. In each of these cases the
  previous lemma says how the flattenings on the remaining simplices
  may be chosen, and moreover, that this choice is unique if $k=2$. If
  $k=3$ then the conditions from the edges $z_0z_4$ and $z_2z_4$
  determine $p_3= p_1-p_0+q_1-q_0$ and $q_3=q_2-q_1$ if
  $(x_0,\dots,x_4)\in\FT^+$, or the approriate analytic continuation
  in general, and the previous lemma then determines the flattening
  $(x_4;p_4,q_4)$.
\end{proof}

\section{Definition of $\themap$}
Let $G=\PSL(2,\C)$ (with the discrete topology). In this section we
describe the map $\themap\co H_3(G;\Z)\to \eprebloch(\C)$.

We first describe the combinatorial representation of elements of
$H_3(G;\Z)$ that we will use. As we will describe, a special case is to
give a closed oriented triangulated $3$--manifold with a flat
$G$--bundle on it. Any element of $H_3(G;\Z)$ can be represented this
way, but we do not want to restrict to this type of
representation. 

We need to clarify our terminology for simplicial complexes. 
\begin{definition}\label{def:simplicial}
  We use the usual concept of \emph{simplicial complex} $K$ except
  that we do not require that distinct simplices have different vertex
  sets (but we do require that closed simplices embed in $|K|$, ie,
  vertices of a simplex are distinct in $K$). The \emph{(open) star}
  of a $0$--simplex $v$ of $K$ is the union of $\tau-\partial\tau$ over
  simplices that have $v$ as a vertex. It is an open neighborhood of
  $v$ and is the open cone on a simplicial complex $L_v$ called the
  \emph{link of $v$}. Note that $L_v$ immerses, but does not
  necessarily embed in $K$ as the boundary of the star of~$v$.
\end{definition}

We really only need quasi-simplicial complexes --- cell complexes
whose cells are simplices with simplicial attaching maps, but no
requirement that closed simplices embed. We discuss this later; our
more restrictive hypothesis is first needed near the end of section
\ref{sec:proof2}, but is eliminated again by Proposition
\ref{prop:local order}.

\subsection{Representing elements of $H_n(G;\Z)$}
\label{subsec:cycles}

\begin{definition}\label{def:cycle}
  An \emph{ordered $n$--cycle} will be a compact $n$--dimensional
  simplicial complex $K$ such that the complement $|K|-|K^{(n-3)}|$ of
  its $(n-3)$--skeleton is an oriented $n$--manifold, together with an
  ordering of the vertices of each $n$--simplex of $K$ so that these
  orderings agree on common faces.  The ordering orients each
  $n$--simplex $\Delta_i$ of $K$, and this orientation may or may not
  agree with the orientation of $|K|-|K^{(n-3)}|$. Let $\epsilon_i=+1$
  or $-1$ accordingly. Then the $n$--cycle $\sum_i\epsilon_i\Delta_i$
  represents a homology class $[K]\in H_n(|K|;\Z)$ called the
  \emph{fundamental class}.
  
  We will also require that the link $L_v$ of each $0$--simplex is
  connected. This can always be achieved by duplicating $0$--simplices
  if necessary.
\end{definition}

One can represent any
element of $H_n(G;\Z)$ by giving an ordered $n$--cycle $K$ and 
labeling the vertices of each simplex $\Delta$ of $K$ by
elements of $G$ so that:
\begin{itemize}
\item Two $G$--labelings $( g_0,\dots,g_k)$ and $(g'_0,\dots,g'_k)$ of an
  ordered $k$--simplex are considered equivalent if there is a $g\in
  G$ with $gg_i=g'_i$ for each $i$;
\item The $G$--labeling of any face of any simplex $\Delta$ is equivalent
  to the $G$--labeling induced from $\Delta$.
\end{itemize}
We will also require that the labels for the vertices of any
$n$--simplex are distinct (we can do this because $G$ is infinite).

We describe two ways of seeing how any element of
$H_n(G,\Z)$ can be described by such data. 
The first is algebraic, and is taken from
\cite{extblochorig}.

\subsubsection{Algebraic description}
We recall a standard chain complex for homology of $G=\PSL(2,\C)$, the
chain complex of ``homogeneous simplices for $G$.''  We will, however,
diverge from the standard by using only non-degenerate simplices,
ie, simplices with distinct vertices --- we may do this because $G$
is infinite.

Let $C_n(G)$ denote the free
$\Z$--module on all ordered $(n+1)$--tuples $\langle
g_0,\dots,g_n\rangle$ of distinct elements of $G$.  Define
$\delta\co C_n\to C_{n-1}$ by $$\delta\langle
g_0,\dots,g_n\rangle=\sum_{i=0}^n(-1)^i\langle g_0,\dots,\hat
g_i,\dots,g_n\rangle.
$$ Then each $C_n$ is a free $\Z G$--module under
left-multiplication by $G$.  Since $G$ is infinite the sequence
$$\cdots \to C_2\to C_1\to C_0\to Z\to 0$$
is exact, so it is a $\Z G$--free resolution of $\Z$.  Thus the chain
complex
$$\cdots\to C_2\otimes_{\Z G}\Z\to C_1\otimes_{\Z G}\Z\to
C_0\otimes_{\Z G}\Z\to 0$$ 
computes the homology of $G$. Note that $C_n\otimes_{\Z G}\Z$ is the
free $\Z$--module on symbols $\langle g_0\co \ldots\co g_n\rangle$, where the
$g_i$ are distinct elements of $G$ and $\langle g_0\co \ldots\co 
g_n\rangle=\langle g_0'\co \ldots\co g_n'\rangle$ if and only if
there is a $g\in G$ with $gg_i=g_i'$ for $i=0,\ldots,n$ (we call these
\emph{homogeneous $n$--simplices for $G$}).

Thus an element of $\alpha\in H_n(G;\Z)$ is represented by a sum
$$
\sum
\epsilon_i\langle g_0^{(i)}\co \ldots\co g_n^{(i)}\rangle
$$
of homogeneous $n$--simplices for $G$ and their negatives (here each
$\epsilon_i$ is $\pm1$).  The fact that this is a cycle means that the
$(n-1)$--faces of these homogeneous simplices cancel in pairs.  We
choose some specific way of pairing canceling faces and form a
geometric quasi-simplicial complex $K$ by taking a $n$--simplex
$\Delta_i$ for each homogeneous $n$--simplex of the above sum and
gluing together $(n-1)$--faces of these $\Delta_i$ that correspond to
$(n-1)$--faces of the homogeneous simplices that have been paired with
each other.

 As already mentioned, a
quasi-simplicial complex actually suffices for our purposes. But we
can obtain a simplicial complex by replacing $K$ by its barycentric
subdivision if necessary (label the barycenter of each simplex by an
arbitrary element of $G$ that differs from the $G$--labels of the
barycenters of all proper faces of that simplex). A standard argument
shows that this does not change the homology class represented by $K$.
We thus get a representation of the homology class in the promised
form.  

\subsubsection{Topological description}\label{subsub:top desc}

Given a space $X$, any element $\alpha\in H_n(X;\Z)$ can be represented as
$f_*[K]$ for some map $f\co |K|\to X$ of an ordered $n$--cycle to $X$.
If $X=BG$ is a classifying space for $G$ then this map is determined
up to homotopy by the flat $G$--bundle $f^*EG$ over $|K|$, so this flat
bundle determines the homology class $\alpha\in H_n(BG;\Z)=H_n(G;\Z)$

Flat $G$--bundles over $|K|$ are determined by $G$--valued $1$--cocycles on
$K$ up to the coboundary action of $G$--valued $0$--cochains (the
$1$--cocycle is the obstruction to extending over $|K^{(1)}|$ a chosen
section of the flat bundle over $K^{(0)}$). We recall the basic
definitions.

 Let $S_q(K)$ be the set of ordered $q$--simplices
of $K$. A $G$--valued $1$--cocyle on $K$ is a map $\sigma\co S_1(K)\to G$
with the cocyle property:
$$
  \sigma\langle v_0 ,v_2\rangle=\sigma\langle
  v_0,v_1\rangle\sigma\langle v_1,v_2\rangle\qquad\text{for }\langle
  v_0,v_1,v_2\rangle\in S_2(K)\,.$$
For convenience we extend the definition of $\sigma$ to
reverse-ordered simplices by 
$$\sigma\langle v_1,v_0\rangle=\sigma \langle v_0,v_1\rangle^{-1}
\qquad\text{for }\langle v_0,v_1\rangle\in S_1(K)\,.
$$
If $\tau\co S_0(K)\to G$ is a $0$--cochain, then its coboundary
action on $1$--cocycles is to replace $\sigma$ by
$$\langle v_0,v_1\rangle \mapsto \tau(v_0)^{-1}\sigma\langle
v_0,v_1\rangle \tau(v_1)\,.$$
A $G$--valued $1$--cocycle $\sigma$ determines a $G$--labeling of
simplices: label a simplex $\langle v_0,\dots,v_k\rangle$ by $(
1,\sigma\langle v_0,v_1\rangle,\dots,\sigma\langle v_0,v_k\rangle)$.
Conversely, a $G$--labeling of simplices determines a 
$1$--cocycle: assign to a $1$--simplex with label
$(g_1,g_2)$ the element $g_1^{-1}g_2\in G$. These correspondences
are clearly mutually inverse.
When a $1$--cocyle $\sigma$ is changed by the coboundary action of a
$0$--cochain $\tau$, the corresponding $G$--labeling $(g_0,\dots,g_k)$ of
a $k$--simplex $\langle v_0,\dots,v_k\rangle$ is replaced by
$(g_0\tau(v_0),\dots,g_k\tau(v_k))$.

Thus, we get our desired representative as in subsection
\ref{subsec:cycles} for a homology class $\alpha\in H_n(G;\Z)$ by
representing $\alpha$ by a flat $G$ bundle over an ordered $n$--cycle
$K$, representing that by a $G$--valued $1$--cycle on $K$, and then
taking the corresponding $G$--labeling of $K$.  We want labels of the
vertices of any $n$--simplex to be distinct, which means that the
$1$--cocycle should never take the value $1$.  Since $G$ is infinite,
this can be achieved by modifying by a coboundary if necessary.

\subsection{Definition of $\themap\co H_3(G;\Z)\to
\eprebloch(\C)$}

Let $K$ be an ordered $3$--cycle. We call a closed path $\gamma$ in $|K|$
a \emph{normal} path if it meets no $0$-- or $1$--simplices of $K$ and
crosses all $2$--faces that it meets transversally.  When such a path
passes through a $3$--simplex $\Delta_i$, entering and departing at
different faces, there is a unique edge $E$ of the 3--simplex between
these faces. We say the path \emph{passes} this edge $E$.

(In the following it will often be necessary to distinguish between a
$1$--simplex of $K$ and the various edges of $3$--simplices that are
identified with this $1$--simplex. To avoid excess notation, we will
often use the same symbol for a $1$ simplex of $K$ and the edges of
$3$--simplices that are identified with it, but we will refer to ``edges''
or ``$1$--simplices'' to flag which we mean.)

Consider a choice of combinatorial flattening $\bfw_i$ for each
simplex $\Delta_i$. Then for each edge $E$ of a simplex $\Delta_i$ of
$K$ we have a log-parameter $l_E=l_E(\Delta_i,\bfw_i)$ assigned.
Recall that this log-parameter has the form $\log z + s\pi i$ where
$z$ is the cross-ratio parameter associated to the edge $E$ of simplex
$\Delta_i$ and $s$ is some integer.  We call ($s$ mod $2$) the
\emph{parity parameter} at the edge $E$ of $\Delta_i$ and denote it
$\delta_E=\delta_E(\Delta_i,\bfw_i)$.

\begin{definition}\label{parity}
  Suppose $\gamma$ is a normal path in $|K|$. The \emph{parity} along
  $\gamma$ is the sum $$\sum_E\delta_E\text{ modulo }2$$
  of the parity
  parameters of all the edges $E$ that $\gamma$ passes. Moreover, if
  $\gamma$ runs in the star of some fixed $0$--simplex $v$ of $K$, then
  the {\em log-parameter along the path} is the sum $$\sum_E\pm
  \epsilon_{i(E)}l_E,$$
  summed over all edges $E$ that $\gamma$
  passes, where:
\begin{itemize}
\item 
$i(E)$ is the index $i$ of the simplex $\Delta_i$ that the edge
  $E$ belongs to and $\epsilon_{i(E)}$ is the coefficient $\pm1$ that
  encodes whether the ordering of the vertices of $\Delta_{i(E)}$
  matches the orientation of $|K|$ or not.
\item
the extra sign $\pm$ is $+$ or $-$ according
as the edge $E$ is passed in a counterclockwise or clockwise
fashion as viewed from the vertex.
\end{itemize}
If $\gamma$ just encircles a $1$--simplex (so the extra signs are
all $+$ or all $-$) we speak of the parity or log-parameter ``around
the $1$--simplex.''
\end{definition}

We assume we have represented an element $\alpha\in H_3(G;\Z)$ by a
$G$--labeled ordered $3$--cycle $K$ as in subsection
\ref{subsec:cycles}.  So to each $3$--simplex $\Delta_i$ of $K$ is
associated a $4$--tuple $\langle g_0^{(i)}\co \dots\co
g_3^{(i)}\rangle$, defined up to left-multiplication by elements of
$G$, and $\epsilon_i=\pm1$ encodes whether the ordering of the
vertices of $\Delta_i$ matches the orientation of $K$ or not.
In the following definition we identify $G=\PSL(2,\C)$ with the isometry
group $\Isom^+(\H^3)$.
\begin{definition}[Flattening of $K$]\label{def:flattening}
  Choose $z\in \partial\overline\H^3$ so $g_0^{(i)}z$,
  $g_1^{(i)}z$, $g_2^{(i)}z$, $g_3^{(i)}z$ are distinct points of
  $\partial\overline\H^3$ for each $i$ (this excludes finitely many points
  $z\in \partial\overline\H^3$).  We then have an
  ideal hyperbolic simplex shape for each simplex $\Delta_i$ of $K$
  and an associated cross ratio $x_i=[g_0^{(i)}z\col g_1^{(i)}z\col
  g_2^{(i)}z\col g_3^{(i)}z]$.  A \emph{flattening of $K$} consists of
  a choice of combinatorial flattenings
  $\bfw_i=\ell(x_i;p_i,q_i)$ of the simplices of $K$ such that the
  parity along any normal path in $K$ is zero and the log-parameter
  around each edge of $K$ is zero. If log-parameter along any normal
  path in the star of each $0$--simplex of $K$ is zero we call it a
  \emph{strong flattening of $K$}.
\end{definition}
\begin{theorem}\label{flattening exists} For each choice of $z$ as in
  the above definition, a strong flattening of $K$ exists.
\end{theorem}
Recall that our $G$--labelled ordered $3$--cycle $K$ was chosen to
represent an element $\alpha\in H_3(G;\Z)$ ($G=\PSL(2,\C)$). The
following theorem finally gives the definition of the map
$\themap\co H_3(G;\Z)\to\ebloch(\C)$.
\begin{theorem}\label{beta exists}  
  For any flattening of $K$ the element
  $\sum_i\epsilon_{i}[x_i;p_i,q_i]\in\eprebloch(\C)$ only depends on
  the homology class $\alpha\in H_3(\PSL(2,\C);\Z)$.  We denote it
  $\themap(\alpha)$.  Moreover, $\themap(\alpha)\in\ebloch(\C)$ and
  $\themap\co H_3(\PSL(2,\C);\Z)\to\ebloch(\C)$ is a homomorphism.
\end{theorem}
In Remark \ref{rem:choosing z} we point out that instead of choosing a
single $z$ in the definition of flattening, we may choose a different $z$ for
each vertex of $K$, and Theorem \ref{beta exists} remains true. This
is useful in practice (Section \ref{sec:3-manifolds}).

\medskip
We now give a brief overview of the proofs.

If the extended Bloch group $\ebloch(\C)$ is to be isomorphic to
$H_3(G;\Z)$ then it must fit in the same Bloch--Wigner short
exact sequence that was given for $H_3(G;\Z)$ in the
Introduction. This is proved in section \ref{sec:ebloch}, and section
\ref{sec:eebloch} is a digression describing a related group. These
sections are independent of the proofs of Theorems \ref{flattening
  exists} and \ref{beta exists}, which take up Sections
\ref{sec:proof} through \ref{sec:proof3}.

First some basic tools are developed in Sections \ref{sec:developing}
and \ref{consequences}. Section \ref{sec:developing} describes a
developing map $\tilde K\to \overline\H^3$ which is helpful for
visualizing the flattening condition. It is used to prove some
preliminary lemmas.  Section \ref{consequences} proves a consequence
of the five term relations, a general ``cycle relation'' which is used
frequently later.

Section \ref{sec:proof} proves Theorem \ref{flattening exists}, the
existence of a strong flattening of $K$. The argument of this section
is combinatorial and depends heavily on \cite{neumann}.

The combinatorial argument is continued at the beginning of the next
section (Section \ref{sec:proof2}) to show that the resulting element
of $\eprebloch(\C)$ is independent of the choice of strong flattening
and lies in $\ebloch(\C)$. It then remains to show that the element of
$\ebloch(\C)$ only depends on the homology class $\alpha\in H_3(G;\Z)$
and not on the $G$--labeled ordered 3--cycle $K$ used to represent it.
In the rest of Section \ref{sec:proof2} we first show that this
element of $\ebloch(\C)$ is unchanged if the $G$--labeling of $K$ is
changed by a coboundary (so it just depends on the flat $G$--bundle and
not on the $1$--cocycle used to represent it) and next that it is
invariant under alteration of the triangulation of $|K|$ by Pachner
moves. Thus, by the end of this section we know that we have an
element of $\ebloch(\C)$ that only depends on the space $|K|$ and the
flat $G$--bundle over it.

This section uses a strong flattening of $K$ and, moreover, it needs
to assume that $K$ is a simplicial rather than quasi-simplicial
complex and that the vertex orderings of its simplices are induced by
a global ordering of the $0$--simplices of $K$. These assumptions will
be removed one by one in the next section, so Theorem \ref{beta
  exists} eventually applies to flattenings that are not necessarily
strong, and also to quasi-simplicial complexes for which the vertex
orderings of the simplices need not be globally induced.

Section \ref{sec:proof3} completes the proof of Theorem \ref{beta
  exists}. First it is shown that the singularities of $|K|$ can
be ``resolved'' to make $|K|$ into a manifold. Since $H_3(G;\Z)$ can
be represented as the bordism group of $3$--manifolds with flat
$G$--bundles, the proof is completed by showing that the element of
$\ebloch(\C)$ is unchanged by elementary bordisms (ie, surgery).
Flattenings and strong flattenings are the same thing when $K$ is a
manifold, so the requirement of strong flattenings dissolves, while
retriangulation using Pachner moves is used to relax the requirements
on $K$ to require only that it be a quasi-simplicial complex with
vertex-ordered simplices.

Finally, the main Theorem \ref{th:main}, stating that $\lambda\co
H_3(\PSL(2,\C);\Z)\to\ebloch(\C)$ is an isomorphism and that
Cheeger--Chern--Simons class is given by Rogers dilogarithm, is proved
in section \ref{sec:ccs}. At this point it follows easily, using work of
Dupont.

Section \ref{sec:unordered} describes the weaker conclusions obtained
if orderings of simplices are not used, confirming a conjecture in
\cite{neumann}.

The next two sections then apply the results to invariants of
$3$--manifolds, and a very brief final section describes what happens
if $\C$ is replaced by a number-field.

\section{Developing map}\label{sec:developing}
Suppose we have a $G$--labeled $3$--cycle $K$.  Choose a generic point
$z\in \partial\overline\H^3$ as in the definition of flattenings
(Definition \ref{def:flattening}). Then  each simplex $\Delta_i$ of
$K$ corresponds to a non-degenerate ideal simplex $\langle g_0^{(i)}z,
g_1^{(i)}z, g_2^{(i)}z, g_3^{(i)}z\rangle$ in $\overline\H^3$. This
simplex is determined up to isometry, since
$(g_0^{(i)},\dots,g_3^{(i)})$ is well defined up to left
multiplication by elements of $G=\PSL(2,\C)=\Isom^+(\H^3)$. This
associates a geometry as an ideal simplex to each simplex of $K$, and
hence, by lifting, also to each simplex of the universal cover $\tilde
K$.

We would like to use this geometry to construct a \emph{developing
  map} $|\tilde K|\to \overline \H^3$. The next lemma says that this is
possible. Note that the ideal simplex corresponding to a $3$--simplex
of $K$ inherits an orientation from $\H^3$ that will not in general
agree with the orientation of the simplex in the manifold $|K-K^{(0)}|$,
so the developing map may be ``folded'' in that some simplices may be
mapped reversing orientation. In particular, adjacent simplices will
map to the same side of a common $2$--face whenever the ideal simplex
orientation of one of the simplices agrees with its orientation in $|K|$
and for the other simplex it differs.

\begin{lemma}
  There is a map $D\co |\tilde K|\to\overline\H^3$ (unique up to
  isometries of $\H^3$) which maps each simplex of the universal
  cover $\tilde K$ of $K$ isometrically and preserving orientation
  with respect to its ideal simplex shape just described. We call it
  the ``developing map.'' The map is equivariant with respect to an
  action of $\pi_1(K)$ on $\overline\H^3$.
\end{lemma}
\begin{proof} Consider the corresponding $1$--cocyle $\sigma$ 
  (subsection \ref{subsub:top desc}) whose value on an unordered
  $1$--simplex $\langle v_1,v_2\rangle$ with label $(g_1, g_2)$ is
  $\sigma(\langle v_1,v_2\rangle)=g_1^{-1}g_2$. Any edge path $\langle
  v_0,v_1\rangle\langle v_1,v_2\rangle \dots\langle
  v_{k-1},v_k\rangle$ determines an element $$\sigma(\langle
  v_0,v_1\rangle)\sigma(\langle v_1,v_2\rangle) \dots\sigma(\langle
  v_{k-1},v_k\rangle)\in G$$
  and this gives a homomorphism $\phi_K$
  from the groupoid of edge paths to $G$.  In particular, it restricts
  to a group homomorphism $\pi_1(K,v_0)\to G$. We thus get an action
  of $\pi_1(K,v_0)$ on $\overline\H^3$.
  
  Choose a $0$--simplex $\tilde v_0$ as a basepoint in the universal
  cover $\tilde K$ of $K$. We can then $G$--label the vertices of
  $\tilde K$ by labeling vertex $\tilde v$ by the element
  $\phi_K(\gamma)$ where $\gamma$ is the image of any edge path from
  $\tilde v_0$ to $\tilde v$ in $\tilde K$. This gives us a lift of
  the $G$--labeling of $K$ to a labeling of $\tilde K$ in which labels
  are absolute rather than just well defined up to left multiplication
  by $G$. Using this labeling, the desired mapping of $|\tilde K|$ to
  $\overline\H^3$ maps a simplex of $\tilde K$ with labels
  $(g_1,\dots,g_4)$ to the ideal simplex $\langle g_1z,\dots,
  g_4z\rangle$. 
  
  To show the equivariance of this map it helps to describe it
  more explicitly. For any two $0$--simplices $w_1,w_2\in \tilde
  K$, let $\phi[w_1,w_2]$ be $\phi_K$ of the image in $K$ of a
  simplicial path $[w_1,w_2]$ from $w_1$ to $w_2$ in
  $\tilde K$.  $\phi[w_1,w_2]$ only depends on $w_1$ and
  $w_2$ since any two paths from $w_1$ to $w_2$ are homotopic.  With
  this notation the
  developing map is  
  $$D(w)=\phi[\tilde v_0,w]z \quad\text{for any }w\in\tilde K^{(0)}.$$
    The action of $\pi_1(K,v_0)$ by covering transformations is such
  that the lift of a closed path representing $\gamma\in \pi_1(K,v_0)$
  is a path $[\tilde v_0,\gamma \tilde v_0]$.

Thus for any $w\in\tilde K^{(0)}$ and $\gamma\in \pi_1(K,v_0)$
we have
$$\begin{aligned}
\phi_K(\gamma)D(w)&=\phi_K(\gamma)\phi[\tilde v_0,w]z=\phi[\tilde
v_0,\gamma\tilde v_0]\phi[v_0,w]z=\\
=&\phi[\tilde
v_0,\gamma\tilde v_0]\phi[\gamma v_0,\gamma w]z
=\phi[v_0,\gamma w]z=D(\gamma w)  
\end{aligned}$$
This proves the equivariance.
\end{proof}
\begin{figure}[ht!]
  \centering
\includegraphics[width=2.5in]{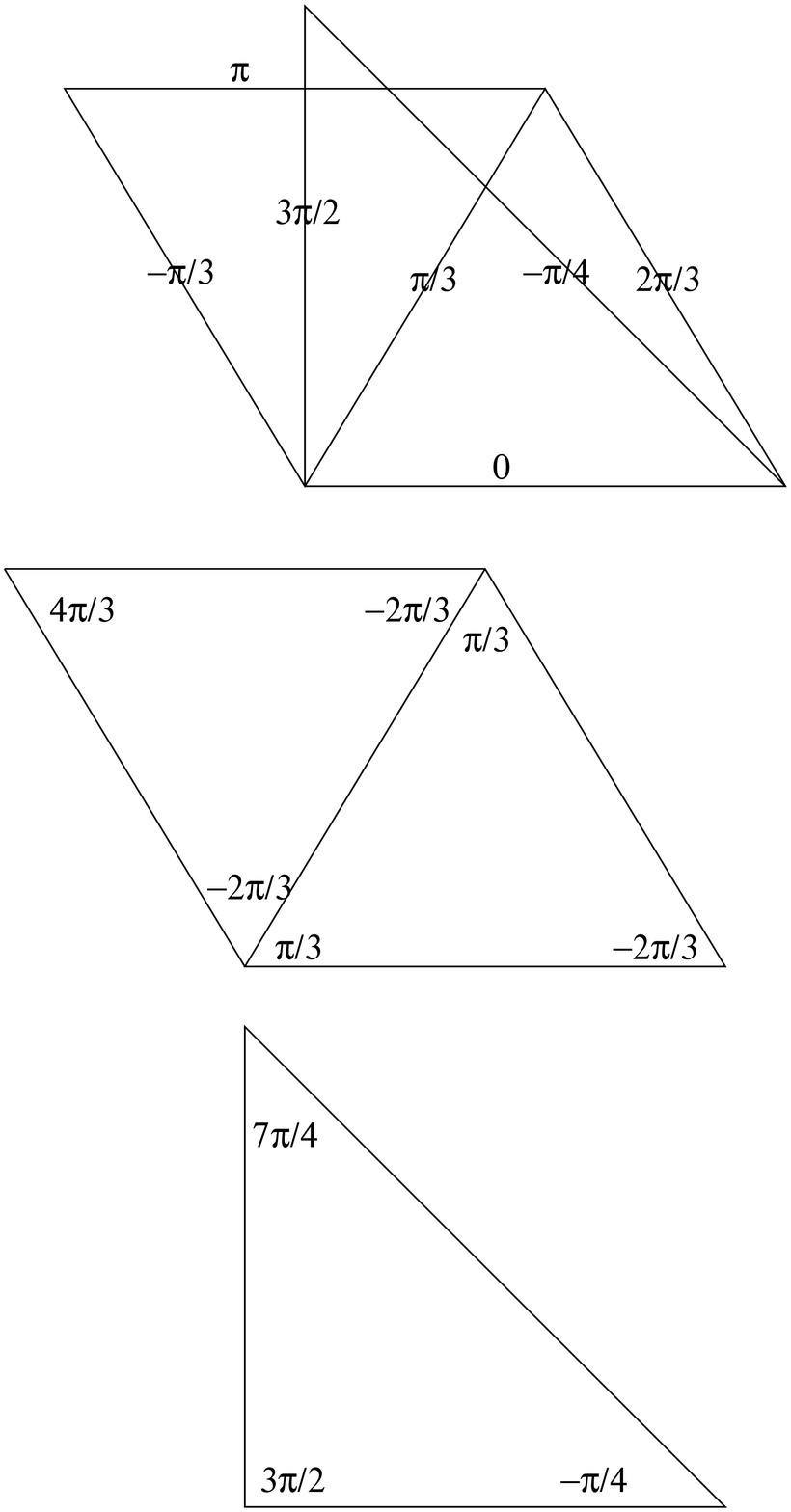}
  \caption{Three triangles in developing image of $N$ in $\C$ showing
    rotation levels of edges and resulting adjusted angles. One triangle is
    from a ``folded'' simplex.}
  \label{fig:flattening}
\end{figure}
The developing map gives a convenient way to visualize what it means
to satisfy the flattening condition in the star neighborhood $N$ of a
$0$--simplex $v$ of $K$. Since $N$ is contractible, it lifts
homeomorphically to the star $\tilde N$ of a lift $\tilde v$ of $v$ in
$\tilde K$. Use the upper half-space model of $\H^3$ and choose the
developing map $D\co |\tilde K|\to \overline \H^3$ so $D(\tilde
v)=\infty$.  Thus for each $3$--simplex of the closure of $\tilde N$,
the other three vertices of the simplex map to points in $\C\subset
\C\cup\{\infty\}=\partial\overline\H^3$, and thus determine a triangle
in $\C$.  Each $2$--simplex incident to $v$ determines a $1$--simplex in
$L_v$ and hence a line segment in $\C$.  For each such segment $S$,
the angle of this line segment from horizontal is well defined in
$\R/\pi\Z$; choose a specific lift $a(S)$ of this number in $\R$ and
call it the \emph{rotation level} for $S$.  Then differences of
rotation levels for segments corresponding to adjacent $2$--faces of a
$3$--simplex determine a flattening in the star neighborhood of $N$
(see Figure~\ref{fig:flattening}).  Specifically:

Suppose we have two adjacent $2$--faces incident to $v$ of a
$3$--simplex $\Delta_i$ of $N$.  Let $S_1$ and $S_2$ be the
corresponding segments in the plane $\C$, taken in the order specified
by the orientation $\Delta_i$ inherits from $K$ (so $S_1$ and $S_2$
are in clockwise order around the triangle in $\C$ determined by
$\Delta$ if the developing map $D$ preserves orientation of $\Delta$
and anticlockwise order otherwise).  Specify the adjusted angle
between these $2$--faces of $\Delta_i$ to be
$\epsilon_i(a(S_2)-a(S_1))$ (where, as usual, $\epsilon_i=\pm1$ is
determined by vertex-order of $\Delta_i$).

This clearly specifies a flattening for the $3$--simplices of $N$ and
the flattening conditions are satisfied at $v$. Conversely, such a
flattening for $N$ plus a rotation level $a(S)$ for one of the line
segments determines $a(S)$ for every segment, so the flattening
determines the rotation level function $a$ up to a multiple of $\pi$.

\smallskip
We will use the developing map to prove a preliminary proposition
towards Theorem \ref{flattening exists}.  That theorem says a strong
flattening of $K$ exists. Since changing flattenings changes log
parameters by multiples of $\pi i$, it must be true that if we make an
\emph{arbitrary} choice of flattenings of the individual ideal
simplices then the sum of log parameters along any normal path in the
star of a $0$--simplex is a multiple of $\pi i$.  We will show, in
fact:
\begin{proposition}\label{prop:preflat}
  Notation as in Definition \ref{def:flattening}.  If we use
  flattenings $\ell(x_i;0,0)$ for the simplices of $K$ then the log
  parameter along any normal path in the star of a $0$--simplex is a
  multiple of $2\pi i$.
\end{proposition}
The proof of this will involve the following related proposition,
which we prove first.
\begin{proposition}\label{prop:preflatp}
  Notation as in Definition \ref{def:flattening}.  If we use flattenings
  $\ell(x_i;0,0)$ for the simplices of $K$ then the parity 
  along any normal path in $K$ is zero.
\end{proposition}
\begin{proof}
  As we follow a normal path the contribution to the parity as we pass
  an edge of a simplex is $0$ if the edge is the $01$, $03$, $12$, or
  $23$ edge of the simplex and the contribution is $\pm1$ if it is the
  $02$ or $13$ edge.  Consider the orientations of the triangular
  faces we cross as we follow the path, where the orientation is the
  one induced by the ordering of its vertices. As we pass a $02$ or
  $13$ edge this orientation changes while for the other edges it does
  not.  Since $K$ is oriented, we must have an even number of
  orientation changes as we traverse the normal path, proving the
  claim.
\end{proof}
\begin{proof}[Proof of Proposition \ref{prop:preflat}]
  Consider a normal path $\gamma$ in the star $N$ of a $0$--simplex $v$
  of $K$. As above, use the upper half-space model of $\H^3$ and
  choose the developing map $D\co |\tilde K|\to \overline\H^3$ so
  $D(\tilde v)=\infty$, where $\tilde v$ is a lift of $v$.
  
  Denote the lift of $\gamma$ to $\tilde N$ by $\tilde\gamma$. The
  path $D\circ\tilde\gamma$ has a ``shadow'' in $\C$. As $\gamma$
  passes an edge of a $3$--simplex of $\tilde N$, the shadow passes the
  corresponding vertex of the corresponding triangle in $\C$, entering
  at one edge of the triangle and exiting at another. The cross-ratio
  parameter corresponding to the $3$--simplex edge (or its inverse,
  depending on vertex ordering and how we are passing the edge) is the
  derivative of the linear map of $\C$ that takes the entering edge of
  the triangle to the exiting edge, fixing the common vertex.  Since
  the path $D\circ\tilde\gamma$ ends up where it started, the product
  of these cross-ratio parameters (each raised to an appropriate power
  $\pm 1$) is $1$, so the appropriately signed sum of their logs is a
  multiple of $2\pi i$.  This sum is precisely the log-parameter along
  the path, except for the adjustments by multiples of $\pi i$ that
  are involved in forming the log-parameters from the logarithms of
  cross-ratio parameters. These adjustments add up to an even multiple
  of $\pi i$ since the parity along the path $\gamma$ to be zero by
  Proposition \ref{prop:preflatp}.
\end{proof}
The following corollary of Propositions \ref{prop:preflat} and
\ref{prop:preflatp} is useful. It is immediate from the definition
of parity.
\begin{corollary}\label{cor:logtoparity}
  For a normal path in the star of a vertex, if the flattening
  condition for log-parameters is satisfied, then so is the parity
  condition.\qed
\end{corollary}

\section{The cycle relation}
\label{consequences}
In this section we prove a general relation called the ``cycle
relation,'' that holds in $\eprebloch(\C)$. We will use it repeatedly
later. It is, in fact, a consequence of the five-term relation
(\ref{five term}) alone and does not involve the transfer relation
(\ref{transfer}).

Let $K$ be a simplicial complex obtained by gluing $3$--simplices
$\Delta_1,\dots,\Delta_n$ together in sequence around a common
\edge{} $E$.  
\begin{figure}[ht!]
\adjustrelabel <-3pt,0pt> {E}{\small$E$}  
\centering
  \includegraphics[width=.25\hsize]{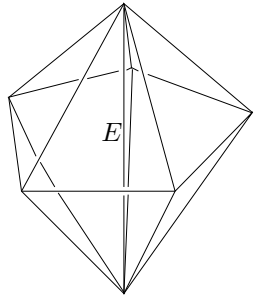}
  \caption{Configuration for the cycle relation}
  \label{fig:cycle}
\end{figure}
Thus, for each index $j$ modulo $n$, $\Delta_j$ is
glued to each of $\Delta_{j-1}$ and $\Delta_{j+1}$ along one of the
two faces of $\Delta_j$ incident to $E$.  Suppose, moreover, that the
vertices of each $\Delta_j$ are ordered such that orderings agree on
the common $2$--faces of adjacent $3$--simplices.

There is then a sequence $\epsilon_1=\pm1$, $\ldots$,
$\epsilon_n=\pm1$ such that the $2$--faces used for gluing all cancel in
the boundary of the $3$--chain $\sum_{j=1}^n\epsilon_j\Delta_j$. (Proof:
choose $\epsilon_1=1$ and then for $i=2,\dots,n$ choose $\epsilon_i$
so the common face of $\Delta_{i-1}$ and $\Delta_i$ cancels.  The
common face of $\Delta_n$ and $\Delta_1$ must then cancel since
otherwise that face occurs with coefficient $\pm2$ in
$\partial\sum_{j=1}^n\epsilon_j\Delta_j$, and $E$ occurs with
coefficient $\pm2$ in $\partial\partial\sum_{j=1}^n\epsilon_j\Delta_j$.)

Suppose now further that a combinatorial flattening $\bfw_i$ has been
chosen for each $\Delta_j$ such that the ``signed sum'' of log
parameters around the \edge{} $E$ vanishes and the same for parity
parameters:
\begin{equation}\sum_{j=1}^n\epsilon_j l_E(\Delta_j,\bfw_j)=0,\quad
\sum_{j=1}^n\epsilon_j \delta_E(\Delta_j,\bfw_j)=0 .
\label{around E}\end{equation}
We think of the \edge{} $E$ as being vertical, so that we can label
the two \edges{} other than $E$ of the common triangle of $\Delta_j$
and $\Delta_{j+1}$ as $T_j$ and $B_j$ (for ``top'' and ``bottom'').
Let $\bfw_j'$ be the flattening obtained from $\bfw_j$ by adding
$\epsilon_j\pi i$ to the log parameter at $T_j$ and its opposite edge
in $\Delta_j$ and subtracting $\epsilon_j\pi i$ from the log parameter
at $B_j$ and its opposite edge in $\Delta_j$.  If we do this for each
$j$ then the total log parameter (and, by Corollary
\ref{cor:logtoparity}, also parity parameter) around any \edge{} of
the complex $K$ is not changed (we sum log-parameters with the
appropriate sign $\epsilon_j$): --- at $E$ no log-parameter has
changed while at every other \edge{} $\pi i$ has been added at one of
the two adjacent $3$--simplices and subtracted at the other.

\begin{lemma}[Cycle relation about $E$]\label{cycle}
With the above notation,
$$\sum_{j=1}^n\epsilon_j[\bfw_j]=
\sum_{j=1}^n\epsilon_j[\bfw_j']~~\in\eprebloch(\C),$$
where we are using $[\bfw]$ as a shorthand for $[\ell^{-1}\bfw]$ 
(ie, $[\bfw]$ means $[z;p,q]$ where $\ell(z;p,q)=\bfw$; see
Lemma \ref{Cover is flattenings}). \end{lemma}

\begin{proof}
  Each of the $3$--simplices $\Delta_i$ has an associated ideal hyperbolic
  structure compatible with the combinatorial flattenings $\bfw_j$.
  This ideal hyperbolic structure is also compatible with the
  flattening $\bfw_j'$.  Choose a realization of $\Delta_1$ as an
  ideal simplex in $\overline\H^3$.  We think of this as a mapping of
  $\Delta_1$ to $\overline\H^3$.  We can extend this to a mapping of
  $K$ to $\overline \H^3$ which maps each $\Delta_j$ to an ideal
  simplex with shape appropriate to its combinatorial flattening.
  Adjacent simplices will map to the same side of their common face in
  $\overline\H^3$ if either their orientations or the signs
  $\epsilon_j$ do not match and will be on opposite sides otherwise.
  The fact that the signed sums of log and parity parameters around the
  \edge{} $E$ are zero guarantees that the identifications match
  up as we go once around $E$.

Note that $K$ has $n+2$ vertices.  We first consider the special case
that $n=3$ and there is an ordering $v_0,\dots,v_4$ of the five
vertices of $K$ that restricts to the given vertex ordering for each
simplex.  We also assume the five vertices of $K$ map to distinct
points $z_0,\dots,z_4$ of $\partial \H^3$.

Each $3$--simplex $\Delta_j$ for $j=1,2,3$ has vertices obtained by
omitting one of the five vertices $v_0,\dots,v_4$.  Denote by $\Delta_4$
and $\Delta_5$ the $3$--simplices obtained by omitting each of the other
two vertices. The fact that the common $2$--faces of the $\Delta_j$
cancel when taking boundary of the chain
$\epsilon_1\Delta_1+\epsilon_2\Delta_2+\epsilon_3\Delta_3$ means that,
up to sign this sum corresponds to three summands of the chain
$\partial\langle v_0,\dots,v_4\rangle =\sum(-1)^i\langle
v_o,\dots,\hat{v_i},\dots,v_4\rangle$. Choose $\epsilon_4$ and
$\epsilon_5$ so that $\sum_{j=1}^5\epsilon_j\Delta_j$ is
$\pm\partial\langle v_0,\dots,v_4\rangle$.  

By Lemma \ref{lem:extend}, we can choose unique combinatorial flattenings
$\bfw_4$ and $\bfw_5$ of $\Delta_4$ and $\Delta_5$ so that the signed
sum of log parameters and parity parameters around any \edge{} of
$K\cup\Delta_4\cup\Delta_5$ is zero.
Note that $\bfw_4$ and $\bfw_5$ do not change if we replace
$\bfw_1,\dots,\bfw_3$ by $\bfw_1',\dots,\bfw_3'$.  By Lemma
\ref{geometric five term} we then have
$$\begin{aligned}
\epsilon_1\bfw_1+\epsilon_2\bfw_2+\epsilon_3\bfw_3 &=
-(\epsilon_4\bfw_4+\epsilon_5\bfw_5)\\
\epsilon_1\bfw_1'+\epsilon_2\bfw_2'+\epsilon_3\bfw_3' &=
-(\epsilon_4\bfw_4+\epsilon_5\bfw_5),
\end{aligned}$$
proving this case.

We next consider the case that for some index $j$ modulo $n$ the
images of $\Delta_j$ and $\Delta_{j+1}$ in $\overline\H^3$ do not
coincide, so their union has five distinct vertices.  By cycling our
indices we may assume $j=1$. Since the orderings of the vertices of
$\Delta_1$ and of $\Delta_{2}$ agree on the three vertices they have
in common, there is an ordering of all five vertices compatible with
both $\Delta_1$ and $\Delta_{2}$. Let $\Delta_0$ be the simplex
determined by the common \edge{} $E$ and the two vertices that $\Delta_1$
and $\Delta_{2}$ do not have in common.  Then there is an
$\epsilon_0=\pm1$ such that the common faces of $\Delta_0$, $\Delta_1$,
and $\Delta_{2}$ cancel in the boundary of the chain
$\epsilon_0\Delta_0+\epsilon_1\Delta_1+\epsilon_{2}\Delta_{2}$.
Choose a flattening $\bfw_0$ of $\Delta_0$ such that
$\epsilon_0l_E(\Delta_0,\bfw_0)+\epsilon_1l_E(\Delta_1,\bfw_1)+
\epsilon_{2}l_E(\Delta_{2},\bfw_{2}) =0$.  Then the relation of the
lemma has already been proved for $\bfw_0$, $\bfw_1$, $\bfw_2$, and
by subtracting this relation from the relation to be proved for
$\bfw_1,\dots,\bfw_n$ we obtain a case of the lemma with one fewer
simplices.  Thus, if we assume the lemma proved for $n-1$ simplices
then this case is also proved.

The above induction argument fails only for the case that there
are $2m$ simplices that alternately ``fold back on each other'' so
that their images in $\overline\H^3$ all have the same four
vertices. The above induction eventually reduces us to this case
(usually with $m=1$).  We must therefore deal with this situation to
complete the proof.  We first consider the case that $m=1$ so $n=2$.
We then have four vertices $z_0,\ldots,z_3$ in
$\partial\overline\H^3$.  We assume the \edge{} $E$ is $z_0z_1$.  Then
the ordering of the vertices of the faces $z_0z_1z_2$ and $z_0z_1z_3$
is the same in each of $\Delta_1$ and $\Delta_2$.  Choose a new point
$z_4$ in $\partial\overline\H^3$ distinct from $z_0,\dots,z_3$ and
consider the ordered simplex with vertices $z_0,z_1,z_2$ ordered as
above followed by $z_4$.  Call this $\Delta_3$.  Similarly make
$\Delta_4$ using $z_0,z_1,z_3$ ordered as above followed by $z_4$.
Choose flattenings of $\Delta_3$ and $\Delta_4$ so that the signed sum
of log parameters for $\Delta_1,\Delta_3,\Delta_4$ around $E$ is
zero.  Then we obtain a three simplex relation of the type already
proved for $\Delta_1,\Delta_3,\Delta_4$ and another for
$\Delta_2,\Delta_3,\Delta_4$, and the difference of these two
relations gives the desired two-simplex relation.

More generally, if we are in the above ``folded'' case with $m>1$ we
can use an instance of the three-simplex relation to replace one of
the $2m$ simplices by two. We then use the induction step to replace
one of these new simplices together with an adjacent old simplex by
one simplex and then repeat for the other new simplex.  In this way we
reduce to a relation involving $2m-1$ simplices, completing the proof.
\end{proof}

\section{Computation of $\ebloch(\C)$}\label{sec:ebloch}
In this section we work out the relationship of the extended Bloch
group with the usual one. The result is Theorem \ref{th:commutdiag}.

We consider the case $n=3$ of the cycle relation (Lemma \ref{cycle}).
Denote the five vertices involved $v_0\dots, v_4$ and assume the three
$3$--simplices meet along the edge $E=v_3v_4$. If we order
$v_0,\dots,v_4$ in this order and give the simplices the inherited
vertex orders then the cycle relation can be written (with the
appropriate relationship among $p_0,p_1,p_2$):
\begin{equation}\begin{aligned}{}
    &[x;p_0,q_0]-[y;p_1,q_1]+[y/x;p_2,q_2]={}\\
    &[x;p_0,q_0-1]-[y;p_1,q_1-1]+
    [y/x;p_2,q_2-1].\label{homo}\end{aligned}\end{equation} This is
true for any choice of $q_0,q_1,q_2$ so long as $p_0,p_1,p_2$ satisfy
the appropriate relation.  Thus if we just change $q_0$ and subtract
the resulting equation from the above we get $$
[x;p_0,q_0]-[x;p_0,q_0']=[x;p_0,q_0-1]-[x;p_0,q_0'-1].  $$
From the
versions of the three-simplex cycle relation with different orderings
of the vertices $v_0,\dots,v_4$ we can similarly derive three versions
of this relation:
\begin{equation}
\begin{aligned}{}
[x;p,q]-[x;p,q']&=[x;p,q-1]-[x;p,q'-1]\\
[x;p,q]-[x;p',q]&=[x;p-1,q]-[x;p'-1,q]\label{three equns}\\
[x;p,q]-[x;p+s,q-s]&=[x;p+1,q-1]-[x;p+s+1,q-s-1]
\end{aligned}
\end{equation}
From these
we obtain:
\begin{lemma}\label{super transfer}
$ [x;p,q]=
pq[x;1,1]-(pq-p)[x;1,0] -(pq-q)[x;0,1]+(pq-p-q+1)[x;0,0]$.
\end{lemma}
\begin{proof}
The first of the relations (\ref{three equns}) implies
$[x;p,q]=[x;p,q-1]+[x;p,1]-[x;p,0]$ and applying this repeatedly shows
\begin{equation}[x;p,q]=q[x;p,1]-(q-1)[x;p,0].\label{dd}\end{equation}
The second equation of (\ref{three
equns}) implies similarly that $[x;p,q]=p[x;1,q]-(p-1)[x;0,q]$ and
using this to expand each of the terms on the right of (\ref{dd})
gives the desired equation.
\end{proof}

Up to this point we have only used consequences of the five-term
relation and not used the transfer relation (\ref{transfer}). We
digress briefly to show that the transfer relation almost follows from
the five term relation.
\begin{proposition}\label{kappa}
If $\eprebloch'(\C)$ and $\ebloch'(\C)$ are defined like
$\eprebloch(\C)$ and $\ebloch(\C)$ but without the transfer relation,
then in $\eprebloch'(\C)$ the element
$\kappa:=[x;1,1]+[x;0,0]-[x;1,0]-[x;0,1]$ is independent of $x$ and
has order $2$. Moreover, $\eprebloch'(\C)={\eprebloch}(\C)\times C_2$
and $\ebloch'(\C)={\ebloch}(\C)\times C_2$, where $C_2$ is the cyclic
group of order 2 generated by $\kappa$.
\end{proposition}
\begin{proof}
  If we subtract equation (\ref{homo}) with $p_0=p_1=q_0=q_1=q_2=1$
  from the same equation with $p_0=p_1=0$, $q_0=q_1=q_2=1$ we obtain
  $[x;1,1]-[y;1,1]-[x;0,1]+[y;0,1]=[x;1,0]-[y;1,0]-[x;0,0]+[y;0,0]$,
  which rearranges to show that $\kappa$ is independent of $x$.  The
  last of the equations (\ref{three equns}) with $p=q=0$ and $s=-1$
  gives $2[x;0,0]=[x;1,-1]+[x;-1,1]$ and expanding the right side of
  this using Lemma \ref{super transfer} gives
  $2[x;0,0]=-2[x;1,1]+2[x;1,0]+2[x;0,1]$, showing that $\kappa$ has
  order dividing 2.

To show $\kappa$ has order exactly 2 we note that there is a
homomorphism $\epsilon\co\eprebloch'(\C)$\break $\to \Z/2$ defined on
generators by $[z;p,q]\mapsto (pq\text{ mod }2)$.  Indeed, it is easy
to check that this vanishes on the lifted five-term relation, and is
thus well defined on $\eprebloch'(\C)$.  Since $\epsilon(\kappa)=1$ we
see $\kappa$ is non-trivial.  Finally, Lemma \ref{super transfer}
implies that the effect of the transfer relation is simply to kill the
element $\kappa$, so the final sentence of the proposition follows.
\end{proof}

\begin{lemma}\label{1-x}
For any $[x;p,q]\in\eprebloch(\C)$ one has $[x;p,q]+[1-x;-q,-p]=2[1/2;0,0]$.
\end{lemma}
\begin{proof}
Assume first that $0<y<x<1$. Then, as remarked in the proof of
Proposition \ref{R}, 
\def\sst{\scriptstyle}
\def\sst{}
$$\begin{aligned}{}
[x;p_0,q_0]&-[y,p_1,q_1]+[{\sst\frac
yx};p_1-p_0,q_2]-[{\sst\frac{1-x^{-1}}{1-y^{-1}}};p_1-p_0+q_1-q_0,q_2-q_1]+{}\\
&[\sst{\frac{1-x}{1-y}};q_1-q_0,q_2-q_1-p_0]=0\end{aligned}$$
is an instance of the lifted five term relation.  Replacing $y$ by
$1-x$, $x$ by $1-y$, $p_0$ by $-q_1$, $q_1$ by $-q_0$, $q_0$ by
$-p_1$, $q_1$ by $-p_0$, and $q_2$ by $q_2-q_1-p_0$ replaces this
relation by exactly the same relation except that the first two terms
are replaced by $[1-y;-q_1,-p_1]-[1-x;-q_0,-p_0]$. Thus subtracting
the two relations gives:
$$[x;p_0,q_0]-[y;p_1,q_1]-[1-y;-q_1,-p_1]+[1-x;-q_0,-p_0]=0.
$$
Putting $[y;p_1,q_1]=[1/2;0,0]$ now proves the lemma for $1/2<x<1$.
But since we have shown this as a consequence of the lifted five term
relation, we can analytically continue it over the whole of $\Cover$.
\end{proof}
\begin{proposition}\label{prop:5.5}
The following sequence is exact:
$$0\to\C^*\stackrel\chi\longrightarrow
{\eprebloch}(\C)\to\prebloch(\C)\to 0$$ where
$\eprebloch(\C)\to\prebloch(\C)$ is the natural map and
$\chi(z):= [z;0,1]-[z;0,0]$ for $z\in \C^*$. 
\end{proposition}
\begin{proof}
Denote $\{z,p\}:=[z;p,q]-[z;p,q-1]$ which is independent of $q$ by the
first equation of (\ref{three equns}).  By Lemma \ref{1-x} we have
$[z;p,q]-[z;p-1,q]=-\{1-z,-q\}$.  It follows that elements of the form
$\{z,p\}$ generate $\Ker\bigl(\eprebloch(\C)\to\prebloch(\C)\bigr)$. Computing
$\{z,p\}$ using Lemma \ref{super transfer} and the transfer relation,
one finds $\{z,p\}=\{z,0\}$ which only depends on $z$.  Thus the
elements $\{z,0\}=\chi(z)$ generate
$\Ker\bigl({\eprebloch}(\C)\to\prebloch(\C)\bigr)$.  If we take equation
(\ref{homo}) with even $p_i$ and subtract the same equation with the
$q_i$ reduced by $1$ we get an equation that says that
$\chi\co\C^*\to \Ker\bigl({\eprebloch}(\C)\to\prebloch(\C)\bigr)$ is a
homomorphism.  We have just shown that it is surjective, and it is
injective because $R\circ\chi$ is the map $\C^*\to \C/\pi^2$ defined
by $z\mapsto \frac{\pi i}2\log z$.
\end{proof}

We can now describe the relationship of our extended groups with the
``classical'' ones.
\begin{theorem}\label{th:commutdiag}
There is a commutative diagram with exact rows and
columns
$$
\begin{CD}
&& 0 &&  0 &&  0 \\
&& @VVV @VVV @VVV \\
0 @>>> \mu^* @>>> \C^* @>>> \C^*/\mu^* @>>> 0\\
&& @V\chi|\mu^* VV @V\chi VV @V\beta VV @VVV\\
0 @>>> \ebloch(\C) @>>> \eprebloch(\C) @>\nu>> \C\wedge\C @>>> 
K_2(\C) @>>> 0\\
&& @VVV @VVV @V\epsilon VV @V=VV\\
0 @>>> \bloch(\C) @>>> \prebloch(\C) @>\nu'>> \C^*\wedge\C^* @>>> 
K_2(\C) @>>> 0\\
&& @VVV @VVV @VVV @VVV\\
&& 0 &&  0 &&  0 && 0
\end{CD}
$$
Here $\mu^*$ is the group of roots of unity and the labeled maps
defined as follows:
$$\begin{aligned}
\chi(z)&=[z;0,1]-[z;0,0]\in
\eprebloch(\C);\\
\nu[z;p,q]&=(\log z +
p\pi i)\wedge(-\log(1-z)+ q\pi i);\\
\nu'[z]&=2\bigl(z \wedge (1-z)\bigr);\\
\beta[z]&=\log z \wedge \pi i;\\
\epsilon(w_1\wedge w_2)&= -2(e^{w_1}\wedge e^{w_2});
\end{aligned}$$
and the unlabeled maps are the obvious ones.
\end{theorem}
\begin{proof}
The top horizontal sequence is trivially exact while the other two are
exact at their first two non-trivial groups by definition of $\ebloch$
and $\bloch$.  The bottom row is exact also at its other two places by
Milnor's definition of $K_2$.  The exactness of the third vertical
sequence is elementary and the second one has just been proved.  The
commutativity of all but the top left square is elementary. A diagram
chase confirms that $\chi$ maps $\mu^*$ to $\ebloch(\C)$ and that the
left vertical sequence is also exact. Another confirms exactness of
the middle row.
\end{proof}

\section{The more extended Bloch group}\label{sec:eebloch}

This section is a digression. We describe a slightly more natural
looking variant of the extended Bloch group, based on the universal
abelian cover $X_{00}$ of $\C-\{0,1\}$. It turns out to be a $\Z/2$
extension of $\ebloch(\C)$. We are not sure of its significance, so we
describe it briefly.

Recall that $\Cover$ consists of four
components $X_{00}$, $X_{01}$, $X_{10}$, and $X_{11}$, of which
$X_{00}$ is naturally the universal abelian cover of $\C-\{0,1\}$.

Let $\liftFT_{00}$ be $\liftFT\cap (X_{00})^5$, so
$\liftFT_{00}=\liftFT_0+2V$ in the notation of Definition
\ref{liftFT}.  As mentioned earlier, $\liftFT_{00}$ is, in fact, equal
to $\liftFT_0$, but we do not need this.

Define $\eeprebloch(\C)$ to be the free $\Z$--module on $X_{00}$
factored by all instances of the relation (we do not need a ``transfer
relation''):
\begin{equation}
\sum_{i=0}^4(-1)^i(x_i;2p_i,2q_i)=0\quad\text{with
$\bigl((x_0;2p_0,2q_0),\dots,(x_4;2p_4,2q_4)\bigr)\in\liftFT_{00}$}.
\end{equation}
As before, we have a well-defined map
$$
\nu\co\eeprebloch(\C)\to\C\wedge\C, 
$$
given by $\nu[z;2p,2q]=(\log z+2p\pi i)\wedge(-\log(1-z)+2q\pi i)$,
and we define
\begin{equation}
\eebloch(\C):=\Ker\nu.
\end{equation}
The proof of Proposition \ref{R} shows:
\begin{proposition}
  The function $R(z;2p,2q):=\calR(z)+\pi i (p\log(1-z)+q\log
  z)-\frac{\pi^2}6$ gives a well defined map $X_{00}\to\C/2\pi^2\Z$
  and induces a homomorphism $R\co \eeprebloch(\C)\to\C/2\pi^2\Z$.
  \qed
\end{proposition}

We can repeat the computations in Section~\ref{consequences}
word-for-word, replacing anything of the form $[x;p,q]$ by
$[x;2p,2q]$, to show that
$\Ker\bigl(\eeprebloch(\C)\to\prebloch(\C)\bigr)$ is generated by
elements of the form $\hat\chi(z):=[z;0,2]-[z;0,0]$.  We thus get:
\begin{proposition} 
The following sequence is exact:
$$
0\to\C^*\stackrel{\hat\chi}\longrightarrow
\eeprebloch(\C)\to\prebloch(\C)\to 0
$$ where
$\eprebloch(\C)\to\prebloch(\C)$ is the natural map and
$\hat\chi(z):= [z;0,2]-[z;0,0]$ for $z\in \C^*$. 
\end{proposition}
\begin{proof}
The only thing to prove is the injectivity of $\hat\chi$ which follows
by noting that $R\bigl(\hat\chi(z)\bigr)=\pi i\log z\in\C/2\pi^2$.\end{proof}
\begin{corollary}
We have a commutative diagram with exact rows and columns:
$$
\begin{CD}
&& 0 &&  0  \\
&& @VVV @VVV  \\
&& \Z/2 @>=>> \Z/2 \\
&& @VVV @VVV  \\
0 @>>> \C^* @>\hat\chi>> \eeprebloch(\C) @>>> \prebloch(\C) @>>> 0\\
&& @VV z\mapsto z^2V @VVV @VV=V \\
0 @>>> \C^* @>\chi>> \eprebloch(\C) @>>> \prebloch(\C) @>>> 0\\
&& @VVV @VVV \\
&& 0 &&  0 && 
\end{CD}
$$
\vbox{and analogously for the Bloch group:
$$
\begin{CD}
&& 0 &&  0  \\
&& @VVV @VVV  \\
&& \Z/2 @>=>> \Z/2 \\
&& @VVV @VVV  \\
0 @>>> \mu^* @>\hat\chi>> \eebloch(\C) @>>> \bloch(\C) @>>> 0\\
&& @VV z\mapsto z^2V @VVV @VV=V \\
0 @>>> \mu^* @>\chi>> \ebloch(\C) @>>> \bloch(\C) @>>> 0\\
&& @VVV @VVV \\
&& 0 &&  0 && \rlap{\phantom{0\hspace{1.3in}0}\qed}
\end{CD}
$$}
\end{corollary}

{\bf Question}\qua Is $\eebloch(\C)$ related to $H_3(\SL(2,\C);\Z)$? Do
$\SL(2,\C)$--labeled ordered $3$--cycles have flattenings using only
$X_{00}$?
\section{Proof of Theorem \ref{flattening exists}}\label{sec:proof}

In this section we will prove Theorem \ref{flattening exists}, which
says that a strong flattening of any $G$--labeled ordered 3--cycle $K$
exists.   

This proof, and the first two lemmas in the next section, depend
heavily on \cite{neumann}. We must first therefore recall some
notation and results from there.  (The results of \cite{neumann} do
not require as strong conditions on $K$ as we have here: there, the
ordering is not needed and $K$ is only required to be a
quasi-simplicial complex.)

Recall that $K$ is an ordered $3$--cycle: the vertices of each
$3$--simplex are ordered with orderings agreeing on common faces. 
%
The underlying space $|K|$ is an oriented $3$--manifold except possibly
at $0$--simplices, where it is topologically the cone on a connected
oriented surface, $L_v$, the link of $v$ (see Definitions
\ref{def:simplicial} and \ref{def:cycle}).

\begin{figure}[ht!]\small\nocolon
\psfrag {e0}{$e_0$}
\psfrag {e1}{$e_1$}
\adjustrelabel <0pt,-1pt> {e2}{$e_2$}
\psfrag {e3}{$e_3$}
\psfrag {e4}{$e_4$}
\adjustrelabel <0pt,-2pt> {e5}{$e_5$}
\centerline{
\includegraphics[width=.3\hsize]{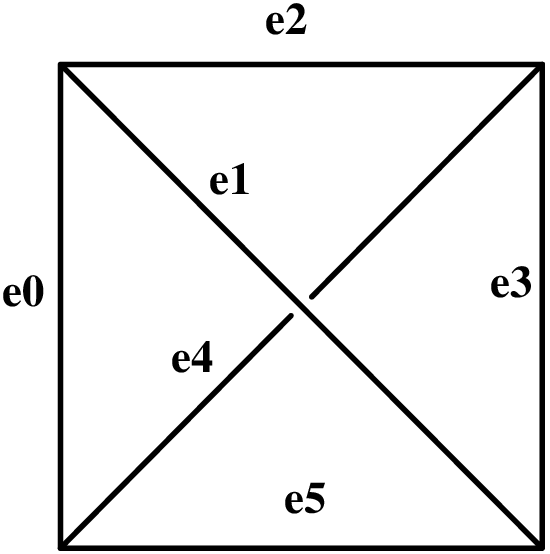}}
\caption{\label{figdelta}}
\end{figure}

To an oriented $3$--simplex $\Delta$ of $K$ we associate a $2$--dimensional 
bilinear space $J_\Delta$ over $\Z$ as follows.  As a $\Z$--module 
$J_\Delta$ is generated by the six edges $e_0,\dots,e_5$ of $\Delta$ 
(see Figure~\ref{figdelta}) with the relations: 
\begin{gather*}
e_i - e_{i+3} = 0 \quad\hbox{for }i=0,1,2.\\
e_0 + e_1 + e_2 = 0.
\end{gather*}

Thus, opposite edges of $\Delta$ represent the same element of 
$J_\Delta$, so $J_\Delta$ has three ``geometric'' generators, 
and the sum of these three generators is zero.  
The bilinear form on $J_\Delta$ is the non-singular skew-symmetric 
form given by
$$
\langle e_0,e_1\rangle = \langle e_1,e_2\rangle = \langle e_2,e_0\rangle = 
-\langle e_1,e_0\rangle = -\langle e_3,e_1\rangle = -\langle e_0,e_2\rangle
= 1. 
$$
Let $J$ be the direct sum $\coprod J_\Delta$, summed over the
oriented $3$--simplices of $K$.  For $i=0,1$ let $C_i$ be the free
$\Z$--module on the \emph{unoriented} $i$--simplices of $K$.  Define homomorphisms
$$
\alpha\co C_0\larrow C_1
\quad\text{and}\quad
\beta\co C_1\larrow J
$$
as follows.
$\alpha$ takes a vertex to the sum of the incident $1$--simplices
(with a $1$--simplex counted twice if both endpoints are at the given vertex).
The $J_\Delta$ component of $\beta$ takes a $1$--simplex $E$ of $K$ to the 
sum of those edges $e_i$ in the edge set $\{e_0,e_1,\dots,e_5\}$ of $\Delta$ 
which are identified with $E$ in $K$.

The natural basis of $C_i$ gives an identification of $C_i$ with its
dual space and the bilinear form on $J$ gives an identification of $J$
with its dual space.  With respect to these identifications, the dual
map
$$
\alpha^*\co C_1\larrow C_0
$$
is easily seen to map a $1$--simplex $E$ of $K$ to the sum of its endpoints,
and the dual map 
$$
\beta^*\co J \larrow C_1
$$
can be described as follows.  To each $3$--simplex $\Delta$ of $K$ we
have a map $j = j_\Delta$ of the edge set $\{e_0,e_1,\dots,e_5\}$ of
$\Delta$ to the set of $1$--simplices of $K$: put $j(e_i)$ equal to the
$1$--simplex that $e_i$ is identified with in $K$.  For $e_i$ in
$J_\Delta$ we have
$$
\beta^*(e_i) = j(e_{i+1}) - j(e_{i+2}) + j(e_{i+4}) - j(e_{i+5})
\quad\hbox{(indices mod 6).}
$$
Let $K_0$ be the result of removing a small open cone neighborhood of
each $0$--simplex $v$ of $K$, so $\partial K_0$ is the disjoint union of the
links $L_v$ of the vertices of $K$.

\begin{theorem}[\cite{neumann}, Theorem 4.2]
The sequence 
$$
\Jscr\co\quad0\larrow C_0\larrow^\alpha C_1\larrow^\beta J
\larrow^{\beta^*} C_1\larrow^{\alpha^*} C_0\larrow 0
$$
is a chain complex. Its homology groups $H_i(\Jscr)$ (indexing the
non-zero groups of $\Jscr$ from left to right with indices $5,4,3,2,1$) are
\begin{gather*}
H_5(\Jscr)=0,\quad H_4(\Jscr)=\Z/2,\quad H_1(\Jscr)=\Z/2,\\
H_3(\Jscr)=\Hscr\oplus H^1(K;\Z/2),\quad
H_2(\Jscr)=H_1(K;\Z/2),
\end{gather*}
where $\Hscr=\Ker(H_1(\partial K_0;\Z)\to H_1(K_0;\Z/2))$.  Moreover,
the isomorphism $H_2(\Jscr)\to H_1(K;\Z/2)$ results by interpreting an
element of $\Ker(\alpha^*)\subset C_1$ as an unoriented 1-cycle in $K$.
\qed 
\end{theorem}

\begin{proof}[Proof of Theorem \ref{flattening exists}: strong
  flattening exists]
Give each simplex $\Delta_i$
of the complex $K$ the flattening $\bfw_i^{(0)}:=\ell(x_i;0,0)$. This
choice will not in general satisfy the conditions for a strong
flattening of $K$, so we need to describe how to modify it so that the
conditions are satisfied.

Recall that to each $3$--simplex $\Delta_i$ of $K$ is associated a sign
$\epsilon_i=\pm1$ that says whether the vertex-ordering of $\Delta_i$
agrees or disagrees with the orientation $\Delta_i$ inherits from the
orientation of $K$.

If $\Delta$ is an ideal simplex and $\bfw=(w_0,w_1,w_2)$ is a
flattening of it, then denote
$$
\xi(\bfw):=w_1e_0-w_0e_1\in J_\Delta\otimes\C.$$
This definition is only apparently unsymmetrical since
$w_1e_0-w_0e_1=w_2e_1-w_1e_2 =w_0e_2-w_2e_0$. 
Denote by $\omega$ the element of $J\otimes\C$ whose
$\Delta_i$--component is $\epsilon_i\xi(\bfw_i^{(0)})$ for each $i$.  That is,
the ${\Delta_i}$--component of $\omega$ is
$-\epsilon_i\bigl(\log(1-x_i)e_0+\log(x_i)e_1\bigr)$.
\begin{lemma}
$\frac1{\pi
i}\beta^*(\omega)$ is an integer class in the kernel of
$\alpha^*$, so it represents an element of the homology group
$H_2(\Jscr)$.  Moreover
 this element in $H_2(\Jscr)$
vanishes, so $\frac1{\pi i}\beta^*(\omega)=\beta^*(\delta)$ for some
$\delta\in J$. 
\end{lemma}
\begin{proof}
Let $\overline J_\Delta$ be defined like $J_\Delta$ but without the
relation $e_0+e_1+e_2=0$, so it is generated by the six edges
$e_0,\dots,e_5$ of $\Delta$ with relations $e_i=e_{i+3}$ for
$i=0,1,2$.  Let $\overline J$ be the direct sum $\coprod\overline
J_\Delta$ over 3--simplices $\Delta$ of $K$.

The map $\beta^*\co J\to C_2$ factors as
$$
\beta^*\co J\larrow^{\beta_1}\overline J\larrow^{\beta_2}C_2,
$$
with $\beta_1$ and $\beta_2$ defined on each component by:
\begin{equation*}
\begin{split}
\beta_1(e_i)&=e_{i+1}-e_{i+2}\\
\beta_2(e_i)&=j(e_i)+j(e_{i+3})
\end{split}
\quad\Biggr\}\text{ for $i=0,1,2$.}
\end{equation*}
Note that $\beta_1(\xi(\bfw))=w_0e_0+w_1e_1+w_2e_2 \in \overline
J_\Delta\otimes\C$.  Thus if $E$ is a $1$--simplex of $K$ then the
$E$--component of $\beta^*(\omega)=\beta_2\beta_1(\omega)$ is the
signed sum of the log-parameters for $E$ in the ideal simplices of $K$
around $E$ and is hence a multiple of $2\pi i$ by Proposition
\ref{prop:preflat}.  That is,
\begin{equation}\frac1{\pi i}\beta^*(\omega)\in 2C_2.\label{even path}
\end{equation} 
The lemma follows, since the isomorphism
$H_2(\Jscr)\to H_1(K;\Z/2)$ is the map which interprets an element of
$\Ker(\alpha^*)$ as an unoriented 1-cycle in $K$, and equation
(\ref{even path}) says this 1-cycle is zero modulo 2.
\end{proof}

Let $\omega':=\omega-\pi i \delta\in J\otimes\C$ with $\delta$ as in
the lemma, so $\beta^*(\omega')=0$.  The $\Delta_i$--component of
$\omega'$ is $\epsilon_i\xi(\bfw_i)$, where $\bfw_i=\ell(x_i;p_i,q_i)$
with the integers $p_i,q_i$ determined by the coefficients occurring
in the element $\delta\in J$. The element $\delta\in J$ is only
determined by the lemma up to elements of $\Ker(\beta^*)$.  We want to
show that for suitable choice of $\delta$, the $\bfw_i$ satisfy the
parity and log-parameter conditions of the definition of strong
flattening (Definition \ref{def:flattening}).  We
will need to review a computation of $H_3(\Jscr)$ from \cite{neumann}.

We define a map $\gamma'\co H_3(\Jscr)\to H^1(\partial
K_0;\Z)=\Hom(H_1(\partial K_0),\Z)$ as follows.  Given elements $a\in
H_3(\Jscr)$ and $c\in H_1(\partial K_0)$ we wish to define
$\gamma'(a)(c)$.  It is enough to do this for a class $c$ which is
represented by a normal path $C$ in the link of some vertex of $K$.
Represent $a$ by an element $A\in J$ with
$\beta^*(A)=0$ and consider the element $\beta_1(A)\in \overline J$.
This element has a coefficient for each edge of each simplex of $K$.
To define $\gamma'(a)(c)$ we consider the coefficients of $\beta_1(A)$
corresponding to edges of simplices that $C$ passes and sum these
using the orientation conventions of Definition \ref{parity}.  It is
easy to see that the result only depends on the homology class of $C$.

We can similarly define a map $\gamma_2'\colon H_3(\Jscr)\!\to\!
H^1(K_0;\Z/2)=\Hom(H_1(K_0),\Z/2)$ by using normal paths in $K_0$ and
taking modulo 2 sum of coefficients of $\beta_1(A)$.

\begin{lemma}[\cite{neumann}, Theorem 5.1]\label{le:ne5.1}
The sequence
$$
0\to H_3(\Jscr)\Mapright{30}{(\gamma',\gamma_2')} H^1(\partial
K_0;\Z)\oplus H^1(K_0;\Z/2)\Mapright{25} {r-i^*}H^1(\partial
K_0;\Z/2)\to 0
$$
is exact, where $r\co H^1(\partial K_0;\Z)\to
H^1(\partial K_0;\Z/2)$ is the coefficient map and the map $i^*\co
H^1(K_0;\Z/2)\to H^1(\partial K_0;\Z/2)$ is induced by the inclusion
$\partial K_0\to K_0$.
\qed\end{lemma}

Returning to the choice of $\delta$ above, assume we have made a
choice so that the resulting flattenings $\bfw_i$ do not lead to zero
log-parameters and parities for normal paths. Taking $\frac1{\pi i}$
times the log-parameters of normal paths leads as above to an element
$c\in H^1(\partial K_0;\Z)$. Similarly, parities of normal paths leads
to an element of $c_2\in H^1(K_0;\Z/2)$. These elements satisfy
$r(c)=i^*(c_2)$.  The lemma thus gives an element of $H_3(\Jscr)$ that
maps to $(c,c_2)$. Subtracting a representative for this element from
$\delta$ gives the desired correction of $\delta$ so the
log-parameters and parities of normal paths with respect to the
corresponding changed $\bfw_i$'s are zero.

This completes the proof that a
strong flattening of $K$ exists: Theorem \ref{flattening exists}. 
\end{proof}

\section{Start of  the proof of Theorem \ref{beta exists}}
\label{sec:proof2}

\begin{lemma}
The choice of strong flattening of $K$ does not affect the
resulting element $\sum_i\epsilon_i [x_i;p_i,q_i] \in\eprebloch(\C)$.
\end{lemma}
\begin{proof}
  If we have a different choice of flattenings $\bfw_i$ satisfying the
  parity and log-parameter conditions for a strong flattening of $K$
  then Lemma \ref{le:ne5.1} implies that the difference between the
  corresponding elements $\delta$ represents $0$ in $H_3(\Jscr)$. It
  is thus in the image of $\beta$.  For $E\in C_2$ the effect of
  changing $\delta$ by $\beta(E)$ is to change the element
  $\sum_i\epsilon_i [x_i;p_i,q_i] \in\eprebloch(\C)$ by the cycle
  relation about $E$ of Lemma \ref{cycle}.  Since this is a
  consequence of the lifted five term relations, the element in
  $\eprebloch(\C)$ is unchanged. 
\end{proof}

Ultimately we will want to see that 
  $\sum_i\epsilon_i [x_i;p_i,q_i]$ is independent of flattening rather
  than strong flattening. This will follow in Section
  \ref{sec:proof3}. The following lemma only needs
  flattening and not strong flattening.
\begin{lemma}\label{le:inbloch}
  Given a flattening of $K$, the element $\sum_i\epsilon_i
  [x_i;p_i,q_i]\in\eprebloch(\C)$ lies in $\ebloch(\C)$.
\end{lemma}
\begin{proof}
For any $\Q$--vector-space $V$, the skew-symmetric bilinear form
$\langle~\rangle$ on $J$ induces a symmetric bilinear map
$$B\co(J\otimes V)\otimes (J\otimes V)\to V\wedge V\,,\quad
(a\otimes v)\otimes (b\otimes w)\mapsto \langle a,b\rangle\, v\wedge
w\,. $$
Theorem 4.1 of \cite{neumann} says that after tensoring with
$\Q$, the bilinear form $\langle,\rangle$ on $\Ker(\beta^*\co
J\otimes\Q\to C_1\otimes\Q)$ induces two times the intersection form
on $H_1(\partial K_0;\Q)=\Ker \beta^*/\Im\beta$. Hence, on
$(\Ker\beta^*/\Im\beta)\otimes V=H_1(\partial K_0)\otimes
V=H_1(\partial K_0;V)$, the above bilinear map $B$ induces the map
$$
B'\co(H_1(\partial K_0)\otimes V)\otimes (H_1(\partial
K_0)\otimes V)\to V\wedge V$$ given by
$$
([a]\otimes v)\otimes ([b]\otimes w)\mapsto 2([a]\cdot [b])\,
v\wedge w$$
where $[a]\cdot[b]$ is intersection form.
  
We take $V=\C$. For our element $w'=w-\pi i \delta \in J\otimes \C$ we
have
$$B(w',w')=2\sum\nolimits_i \epsilon_i(\log x_i + p_i\pi
i)\wedge(-\log(1-x_i)+ q_i\pi i)\in\C\wedge\C\,.$$
Thus, we want to
show that $B(w',w')=0$.  If we have a strong flattening, $w'$ is in
$\Ker(\beta^*)\otimes \C$ and represents zero in $H_1(\partial
K_0;\C)$. Thus $B(w',w')=B'([w'],[w'])=B'(0,0)=0$, as desired.  If the
flattening is not strong, then rather than $w'$ representing zero in
$H_1(\partial K_0;\C)$, it represents $\pi i$ times an integral
homology class $\alpha$ say, and we still have
$B(w',w')=B'([w'],[w'])=(\alpha\cdot\alpha)\pi i\wedge \pi i=0$.
\end{proof}

At this point we have shown that the element $\sum_i\epsilon_i
[x_i;p_i,q_i]$ determined by a strong flattening of $K$ lies in
$\ebloch(\C)$ and is independent of choice of strong flattening. In
the remainder of this section we show that it depends only on $|K|$
and the flat $G$--bundle over it.

\begin{lemma}\label{le:changez}
  Changing the choice of the point $z\in\partial\overline\H^3$ in the
  definition of strong flattening (Definition \ref{def:flattening})
  does not change the element $\sum_i\epsilon_i
  [x_i;p_i,q_i]\in\ebloch(\C)$.
\end{lemma}
\begin{proof}
  This is a special case of the following lemma, since if we change
  the $1$--cocycle corresponding to a $G$--labeling by the coboundary of
  the constant $0$--cochain with value $g$, the effect is the same as
  replacing $z$ by $gz$ (see subsection \ref{subsub:top desc}).
\end{proof}
\begin{lemma}\label{le:changeG}
  If we change the $G$--labeling of $K$ by changing the corresponding
  $1$--cocycle by a coboundary (subsection \ref{subsub:top desc}) then
  the element $\sum_i\epsilon_i [x_i;p_i,q_i]\in\ebloch(\C)$ does not
  change.
\end{lemma}
\begin{proof}
  The corresponding result is known for the element in $\bloch(\C)$,
  so the change is in $\Ker(\ebloch(\C)\to\bloch(\C))=\Q/\Z$. Thus,
  if we knew that the coboundary action of $0$--cochains was continuous
  using the standard topology of $G=\PSL(2,\C)$ the Lemma would
  follow. This continuity seems ``self-evident'' but we do not know an
  easier proof than what follows, which directly proves the local
  constancy of the element of $\ebloch(\C)$ under this action.
  
  It suffices to prove the lemma for the change given by the
  coboundary action of a $0$--cochain $\tau$ that takes the value $1$
  on all $0$--simplices except one. Denote that one $0$--simplex by $v$.
  The effect is that the ideal simplex corresponding to a simplex
  $\Delta$ of $K$ is unchanged if $v$ is not one of its vertices,
  while if it has vertices $v,v_1,v_2,v_3$ with $G$--labels
  $g,g_1,g_2,g_3$ before the change, then the corresponding ideal
  simplex has vertices $gz,g_1z,g_2z,g_3z$ before the change and
  $g\tau(v)z,g_1z,g_2z,g_3z$ after. We assume $\tau(v)$ is
  sufficiently close to $1$ that none of these latter simplices are
  degenerate (two vertices equal).
  
  Since the flattening condition on $K$ 
  is a discrete condition (it can only fail by multiples of $\pi i$),
  it will stay valid if we vary strong flattenings of simplices continuously
  as $\tau(v)\in G$ varies (considering $G=\PSL(2,\C)$ with its usual
  topology, rather than as a discrete group).  Thus we get a
  continuous family of strong flattenings of $K$ as $\tau(v)$ varies in a
  neighborhood of $1\in G$.  We must show that any two of them lead to
  equal elements of $\ebloch(C)$, ie, they can be related using
  lifted five-term relations. Since only simplices in the star $N$ of
  $v$ are affected, it suffices to show that
\begin{equation}\label{eq:change}
\sum_{\Delta_i\subset
    N}\epsilon_i([x_i;p_i,q_i]-[x'_i;p'_i,q'_i])=0\in\eprebloch(\C),
\end{equation}
where $\ell([x_i;p_i,q_i])$ and $\ell([x'_i;p'_i,q'_i])$ are the
flattenings of $\Delta_i$ before and after changing the $G$--label of
$v$ by $\tau(v)$.
  
Denote the $0$--simplices in the simplicial link of $v$ by $v_1,\dots,
v_m$, and let their $G$--labels be $g_1,\dots,g_m$.  For each ordered
$2$--simplex $\langle v_i,v_j,v_k\rangle$ in the simplicial link of $v$
we wish to give a lifted five-term relation based on the five points
$gz$, $g\tau(v)z$, $g_iz$, $g_jz$, $g_kz$ so that when we sum these
relations (with appropriate signs) we get the left side of equation
(\ref{eq:change}).
  
We already have flattenings of the two simplices $$\langle
gz,g_1z,g_2z,g_3z\rangle\quad\text{and}\quad\langle
g\tau(v)z,g_1z,g_2z,g_3z\rangle.$$
We wish to find flattenings of the
other three simplices $$\langle gz,g\tau(v)z,g_jz,g_kz\rangle,\quad
\langle gz,g\tau(v)z,g_iz,g_kz\rangle,\quad \langle
gz,g\tau(v)z,g_iz,g_jz\rangle$$
so that we have an instance of the
five-term relation.
  
Consider the appropriate signed sum of the three adjusted angles about
the $1$--simplex $\langle g_jz,g_kz\rangle$ in the the flattening
condition of the five-term relation. This sum is a multiple of $\pi$
which should be zero. The two adjusted angles at this edge in $\langle
gz,g_iz,g_jz,g_kz\rangle$ and $\langle
g\tau(v)z,g_iz,g_jz,g_kz\rangle$ almost cancel in the sum. Since we
may assume $\tau(v)$ is close to $1$, the angle at the edge $\langle
g_jz,g_kz\rangle$ in the simplex $\langle
gz,g\tau(v)z,g_jz,g_kz\rangle$ will be small. Thus the adjustment to
this angle must be zero.  Similarly for the $\langle g_iz,g_kz\rangle$
edge of $\langle z,\tau(v)z,g_iz,g_kz\rangle$ and the $\langle
g_iz,g_jz\rangle$ edge of $\langle z,\tau(v)z,g_iz,g_jz\rangle$. In
each of these three simplices the adjusted angle $\langle
z,\tau(v)z\rangle$ equals the adjusted angle at the opposite edge and
is hence small, so the sum of these adjusted angles must be $0$ (since
it is a multiple of $\pi$) so the flattening condition holds also at
the edge $\langle z,\tau(v)z\rangle$.
  
It remains to choose the angle adjustments at edges of the form
$\langle \tau(v)z,g_iz\rangle$ of these simplices.  We use the
description of flattenings in terms of rotation levels on edges
described at the end of section \ref{sec:developing}. So we
position the developing map restricted to $N$ so that vertex $v$ goes
to $gz=\infty$. Then $g\tau(v)z$ is close to $\infty$, ie, very far
from the points $g_\nu z$ for $\nu=1,\dots,m$. Thus the segments
joining $g\tau(v)z$ to the points $g_\nu z$ will be almost mutually
parallel, so we can choose the rotation levels of these segments to be
almost equal to each other. With such a choice it is clear by
inspection that the remaining flattening conditions for the desired
five-term relations hold. When we sum these five-term relations with
appropriate signs, each flattened simplex with vertices of the form
$gz=\infty, g\tau(v)z, g_iz,g_jz$ occurs in two of the relations and
cancels and we are left with the desired relation (\ref{eq:change}).
  
This proves the existence of the desired lifted five-term relations
when $\tau(v)$ is close to $1$. This shows that the element of
$\ebloch(\C)$ is locally constant under the coboundary action of
$0$--cochains, and since the space $G^{K^{(0)}}$ of $0$--cochains is
connected, the Lemma follows.
\end{proof}
\begin{remark}\label{rem:choosing z}
  The above lemma can be interpreted to say that, instead of choosing
  a single $z\in\partial\H^3$ in Definition \ref{def:flattening} to
  obtain the ideal simplex shapes for the simplices of $K$, we can
  choose a different $z$ for each $0$--simplex of $K$, and the element
  of $\ebloch(\C)$ is unaffected. This will be important when
  discussing the applications to $3$--manifolds in Section
  \ref{sec:3-manifolds}.
\end{remark}

We have been assuming that $K$ is a simplicial rather than
quasi-simplicial complex, but so far we have not really used this. We
will initially prove Theorem \ref{beta exists} under the following
assumption, which we will eliminate again in Proposition
\ref{prop:local order}.
\begin{assumption*}
  $K$ is a simplicial complex and the vertex orderings of the
  simplices of $K$ are inherited from an ordering of the set of
  $0$--simplices of $K$.
\end{assumption*}
\begin{lemma}\label{le:changeorder}
  Changing the ordering of the $0$--simplices of $K$ does not change
  the element $\sum_i\epsilon_i [x_i;p_i,q_i]\in\ebloch(\C)$.
\end{lemma}
\begin{proof}
  The full permutation group on $K^{(0)}$ is generated by
  transpositions on adjacent elements with respect to the ordering, so
  we need only consider such transpositions. If the two $0$--simplices
  are not joined by a $1$--simplex then the the transposition has no
  effect on the vertex-ordering of any $3$--simplex so there is nothing
  to prove. Assume therefore that $v$ and $w$ are two $0$--simplices
  that are adjacent in the ordering of $K^{(0)}$ and are joined by a
  $1$--simplex. Only the $3$--simplices that have $\langle v,w\rangle $
  as an edge have their vertex-ordering changed by the transposition,
  so we just consider these simplices. The configuration of these
  simplices is as in Lemma \ref{cycle} and the same inductive argument
  used in the proof of that lemma lets us deduce the lemma from the
  case that there are just three $3$--simplices about this $1$--simplex
  and their vertex orderings are induced from an ordering of the five
  $0$--simplices involved. In the proof of Lemma \ref{cycle} it is
  shown that there is then a lifted five-term relation which replaces
  these three simplices by two.  After doing so, $v$ and $w$ no longer
  joined by an edge, so the transposition no longer has an effect on
  the vertex-order of any three-simplex, so the lemma follows.
\end{proof}

Finally for this section, we show that we can change the triangulation
of $K$ by Pachner moves without changing the element of $\ebloch(\C)$.
These are the moves ($1\to 4$, $2\to 3$, $3\to 2$ and $4\to 1$ on
numbers of simplices) that replace a union of $3$--simplices of $K$
that is combinatorially equivalent to part of the boundary of a
$4$--simplex by the union of the complementary set of $3$--simplices in
the boundary of the $4$--simplex.  Only the $1\to 4$ Pachner move adds
a new $0$--simplex. In this case, we insert the new $0$--simplex
anywhere in the ordering of $0$--simplices and, if the old simplex had
$G$--labels $(g_0,g_1,g_2,g_3)$, we give the four new simplices
$G$--labels $(g_1,g_2,g_3,g)$, $(g_0,g_2,g_3,g)$, $(g_0,g_1,g_3,g)$,
and $(g_0,g_1,g_2,g)$, for some $g\ne g_0,g_1,g_2,g_3$.  We are again
using the assumption of a global ordering of the $0$--simplices of $K$:
if we only had local vertex-orderings of simplices, we could not
guarantee that the two $3$--simplices produced by a $3\to 2$ Pachner
move have vertex-orderings compatible with adjacent simplices.  
\begin{lemma}\label{le:pachner}
  A change of triangulation of $K$ by Pachner moves does not change
  the represented element $\sum_i\epsilon_i
  [x_i;p_i,q_i]\in\ebloch(\C)$.
\end{lemma}
\begin{proof}
  We want to leave unchanged the flattenings on the unaltered
  simplices of $K$ and put flattenings on the new simplices to get a
  flattening of the new complex $K'$. We can do this if and only if
  the flattenings of the changed simplices (the ones in $K$ that have
  been replaced and the ones they have been replaced by) give a
  flattening of the boundary of the $4$--simplex (in the sense that log
  and parity parameters around edges are zero). In particular, the
  change in the element of $\ebloch(\C)$ is then an instance of the
  lifted five-term relation, and hence zero. 
  We thus need to know that any flattening defined on part of the
  boundary of a $4$--simplex can be extended over the whole
  $4$--simplex.  This is the content of Lemma \ref{lem:extend}.
\end{proof}
It follows from Pachner \cite{pachner} that any two simplicial
triangulations of $|K|$ are related by Pachner moves, so at this point
we know that the element of $\ebloch(\C)$ is determined just by $|K|$
and the flat $G$--bundle over it. We will actually only need this in
the case that $|K|$ is a manifold, which is explicit in Pachner's
work. In the next section we reduce to the case that $|K|$ is a
manifold and use this to complete the proof of Theorem \ref{beta
  exists}.

\section{Completion of proof of Theorem \ref{beta exists}}\label{sec:proof3}

In this section we show that we can always resolve the singularities
of $|K|$ to replace it by a triangulated $3$--manifold and then discuss
surgery on this $3$--manifold to complete the proof of Theorem
\ref{beta exists}.

\begin{figure}[ht!]
  \centering
  \includegraphics[width=.5\hsize]{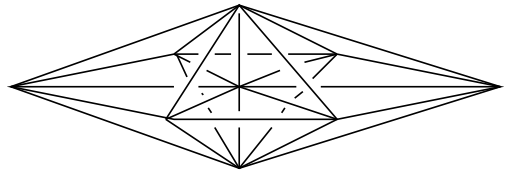}
  \caption{Double cone on a disc}
  \label{fig:cone}
\end{figure}

The underlying space of $K$ is a $3$--manifold except at $0$--simplices,
where it may look locally like the cone on a closed surface of genus
$\ge 1$. Suppose this occurs at some $0$--simplex $v$. Choose a
non-separating simple closed simplicial curve $C$ in the link $L_v$.
This curve is determined by an open disk $D$ in the star of $v$. We
slice $|K|$ open along $D$ and splice in a complex as in
Figure~\ref{fig:cone}, gluing the top and the bottom of that complex to
the two copies of $D$ resulting from slicing $|K|$.
This introduces a new $0$--simplex $v'$ in $K$.  We give $v'$ a
$G$--label distinct from that of each adjacent $0$--simplex.  We order
$0$--simplices by inserting $v'$ anywhere (eg, as a new maximal
element). We must check that we can put flattenings on the inserted
$3$--simplices so that the conditions for a strong flattening of our
complex still hold. 

We first show that it suffices to satisfy the log-parameter flattening
conditions.  Indeed, by Corollary \ref{cor:logtoparity}, if we satisfy
the logarithmic flattening conditions, the parity conditions hold
around edges. Parity along normal paths then gives a homomorphism from
$H_1(|K-K^{(0)}|)$ to $\Z/2$. If we denote the new complex by $K'$ and
the curve in $K'$ given by the two vertical $1$--simplices by $C$, then
$H_1(|K-K^{(0)}|)=H_1(|K'-K^{(0)}|- C)$ which surjects to
$H_1(|K'-(K')^{(0)}|)$ by the inclusion, so the parity condition will
hold for all normal paths.

We satisfy the log-parameter strong flattening condition as follows.  The
added $3$--simplices come in pairs with identical $G$--labelings,
adjacent to each other above and below the central disk in the
picture. We will give the two simplices of a pair identical
flattenings so they cancel in the sum $\sum_i\epsilon_i
[x_i,p_i,q_i]\in\ebloch(\C)$ (they appear with opposite signs in this
sum).  Flatten all the added simplices in the top layer in order to
satisfy the flattening conditions at the top edges in the picture
(there are $n$ degrees of freedom in doing this, where $n$ is the
number of $1$--simplices around the curve $C$). We then need only
verify the flattening condition around the central vertical $1$--simplex.
Since the sum of the three log-parameters of a flattened simplex is
zero, the appropriately signed sum of log parameters around the
central edge is the negative of the signed sum of the $2n$
log-parameters of these simplices at the top edges in the picture.
This, in turn, is equal to the log-parameter along a normal path in
the star of $v$ determined by a path in the link $L_v$ that runs
parallel to $C$ on one side. It is thus zero by the strong flattening
condition for such normal paths.

(An alternate argument is to look at the developing image centered at
$v$ as in section \ref{sec:developing}. The flattenings before
we change $K$ give us rotation levels on the edges of the image in
$\C$. After changing $K$ we have a new point in $\C$ at the image of
$v'$ and $n$ segments incident on it, corresponding to the
$2$--simplices of the central horizontal disk in the picture.  Assign
any rotation levels to these $n$ segments. As in section
\ref{sec:developing}, this determines flattenings of the
$3$--simplices so that the flattening conditions are satisfied.)

We call the above procedure ``blowing up'' at the $0$--simplex $v$.
Note that it reduces the genus of the link of $v$. By repeated blowing
up we can thus make $|K|$ into a $3$--manifold. The following lemma 
improves on this.
\begin{lemma}
  If $K$ is just flattened (rather than strongly flattened) then we
  can blow up to make $|K|$ into a $3$--manifold without changing the
  represented element in $\ebloch{\C}$.  In particular, Theorem
  \ref{beta exists} holds for flattenings if it holds for strong
  flattenings.
\end{lemma}
\begin{proof}
  To apply the above argument to to blow up and simplify the link of a
  vertex we need a separating curve $C$ in the link so that the
  log-parameter along a parallel normal curve in the link $L_v$ is
  zero.  Suppose the log-parameter flattening conditions are satisfied
  around edges, but not necessarily along normal curves.
  Log-parameters along normal curves then give a homomorphism
  $H_1(L_v;\Z)\to \Z\pi i\subset \C$. It is well known that for any
  closed oriented surface $S$ and non-trivial homomorphism $H_1(S;\Z)
  \to \Z$ there exist elements in the kernel represented by
  non-separating simple closed curves (in fact, the mapping class
  group acts transitively on the set of surjective homomorphisms
  $H_1(S;\Z) \to \Z$). Thus, there exists a curve to use for blowing
  up, so the lemma follows.
\end{proof}

At this point we may assume our element in $H_3(G;\Z)$ is represented
by a triangulated $3$--manifold with flat $G$--bundle over it, or
equivalently, a $3$--manifold $|K|$ with a homotopy class of mappings
$|K|\to BG$. Since $H_3(BG;\Z)=\Omega_3(BG)$ (oriented bordism), and
bordism is generated by homotopy (ie, maps of $|K|\times I$) and
handle addition, it suffices to show that when we modify $K$ by
surgery resulting from a handle-addition to $|K|\times I$ the
represented element in $\ebloch(\C)$ is not changed.

It suffices to consider adding $1$--handles and $2$--handles, since any
bordism of connected $3$--manifolds is the composition of bordisms that
add $1$, $2$, and $3$ handles in that order, and the latter are
inverses of bordisms that add $1$--handles.

We first consider adding a $1$--handle. Combinatorially, homotoping
$|K|$ in $BG$ means changing the $G$--valued $1$--cocycle by the
coboundary action. We make such a change so that the handle addition
can be realized by gluing a $4$--simplex $\Delta^4$ to $K$ by gluing
two of the $3$--faces of $\Delta^4$ to disjoint $3$--simplices
$\Delta_1^3$ and $\Delta_2^3$ in $K$. The common $2$--simplex of the
two $3$--simplices in $\partial \Delta^4$ determines $2$--faces of
$\Delta_1^3$ and $\Delta_2^3$ that will be glued together in this
construction. Our initial change of the $G$--valued $1$--cocycle
arranges for the cocycle to match on these two simplices.

The resulting surgery of $K$ is the result of first removing the
interiors of $\Delta_1^3$ and $\Delta_2^3$, then gluing the resulting
boundary components along the $2$--simplices mentioned above to give a
manifold with a single $2$--sphere boundary triangulated as the
suspension of a triangle, and finally gluing in a $3$--ball
triangulated as three $3$--simplices meeting along a common
$1$--simplex. Call the resulting triangulated manifold $K'$.

We must make the construction compatible with orderings and provide
suitable flattenings. By Lemma \ref{le:changeorder} we may
order the vertices of $K$ as $v_1,v_2,\dots$ so that the vertices of
$\Delta_1^3$ are $v_1,v_3,v_5,v_7$ and of $\Delta_2^3$ are
$v_2,v_4,v_6,v_8$ and the $2$--simplices that are identified are given
by the first three vertices of each in the given order. Then the
vertices of $K'$ are $v_1=v_2,v_3=v_4,v_5=v_6,v_7,v_8,\dots$ and the
common $3$--simplices of $K$ and $K'$ have not changed their vertex
orderings.

Each of the two $3$--simplices $\Delta_1^3$ and $\Delta_2^3$ had a
flattening before the surgery, and we use lemma \ref{lem:extend} to
define the flattening on the three new $3$--simplices. The
log-parameter part of the flattening condition is then satisfied in
$K'$, but the parity condition may not be. The surgery has added a
generator to $H_1(K,\Z)$ represented by any normal path that passes
once through the $2$--simplex that we glued. We must check the parity
condition along one such path. If it fails, then before performing the
surgery we change flattenings of $K$ by the cycle relation of Lemma
\ref{cycle} about the $1$--simplex $v_1v_3$. After doing so the parity
part of the flattening condition will be satisfied for $K'$.

Finally, we must consider adding a $2$--handle along a simplicial curve
in $K$. We can realize this by doing the reverse of a blow-up (a
``blow-down'') to collapse the curve to a single vertex $v$ with link
a torus, and then blowing up again using a longitude of the curve we
collapsed. We need to show, therefore, that if $K$ contains a curve
that maps to the trivial element in $G$ under the homomorphism
determined by the $G$--valued $1$--cocycle, then we can (after adjusting
the triangulation) find a flattening that lets us perform a blow-down
to collapse the curve, followed by a blow-up using a longitude of the
curve.

By retriangulating we may assume the curve $C$ in question is length 2
and the $3$--simplices that have a $1$--simplex in common with it are
configured as in Figure~\ref{fig:cone} (with, for concreteness, $6$
simplices around each edge of $C$, as in the figure). In particular,
the bottom and top vertices in the figure are the same $0$--simplex $v$
in $K$, and, since the $1$--cocyle along the curve gives the trivial
element of $G$, if we $G$--label the simplices of the figure all with
the same label for $v'$ then each pair of vertically adjacent
$3$--simplices has the same labeling.  The rest of the star of $v$ will
be the cone on an annulus. In particular, the $1$--simplices radiating
from $v$ on the top surface of the figure are distinct from the
$1$--simplices radiating from $v$ on the bottom surface.

We must verify that we can modify the flattenings so that vertically
adjacent $3$--simplices have the same flattening, in which case we
can perform the desired blow-down.

Denote the $1$--simplices radiating from $v$ in the bottom surface of
the figure $e_1,\dots,e_6$ in order, and the $1$--simplices above these
in the central disk $e'_1,\dots,e'_6$ (with $e'_i$ above $e_i$).
Denote by $\Delta_i$ the $3$--simplex in the bottom layer with edges
$e_i$ and $e_{i+1}$ (indices mod $6$). By applying suitable multiples
of the cycle relation of Lemma \ref{cycle} around $1$--simplices $e_1$
and $e_2$ we can adjust the flattening of $\Delta_1$ arbitrarily, so
we can make it match the $3$--simplex above it. The flattening
condition at $1$--simplex $e'_2$ now implies that the log-parameter of
$\Delta_2$ at edge $e'_2$ matches the simplex above it, and by
applying the cycle relation around $e_3$ we can fix up
the flattening of $\Delta_2$ at edge $e'_3$. Continuing this way, we
fix $\Delta_3, \dots,\Delta_5$, at which point the flattening
conditions about $e'_5$ and $e'_6$ imply that $\Delta_6$ also has its
desired flattening.

We can now do the blow-down by removing the interior of the complex
shown in the figure from $K$ and then gluing the top and bottom
surfaces. The result is a complex $K'$ in which the link of
$0$--simplex $v$ is a torus, and if we can blow up using a
complementary curve in this torus then our handle addition is
complete. However, to blow up we need the log-parameter along the
appropriate normal curve to be zero, which need not be the case.  We
describe how to remedy this.

Suppose we are in the situation just before we do the blow-down, so
the flattening have been matched on pairs of $3$--simplices in the
figure. If we apply the cycle relation about the lower vertical
$1$--simplex joining $v$ to $v'$ we destroy this matching. The
procedure described above to make flattenings match again then applies
the cycle relation once about each of the $1$--simplices
$e_1,\dots,e_6$. We claim that the end result is to change the
log-parameter along the normal curve we are interested in by $2\pi i$.

\begin{figure}[ht!]
  \centering
    \includegraphics[width=.8\hsize]{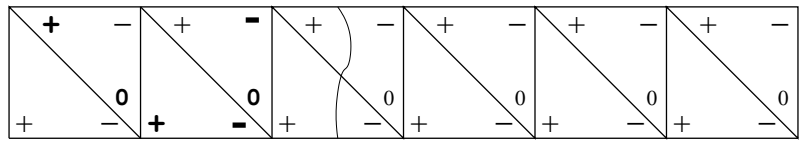}
 \caption{Left and right sides are identified to form an annulus. 
   Signs indicate change of flattening by $\pm\pi i$ and $0$ signifies
   that flattening has been changed by both $\pi i$ and $-\pi i$ so
   total change is zero. The log-parameters affected by applying the
   cycle relation about a single $e_i$ are bold.}
  \label{fig:annulus}
\end{figure}

Indeed, recall that the full star neighborhood of the curve $C$
consists of the complex of Figure~\ref{fig:cone} glued to the cone on
an annulus. In Figure\ \ref{fig:annulus} we show a possible
triangulation of the annulus (we may assume it is triangulated this
way, since we may choose the triangulation). The figure also shows
how log-parameters are changed by the above procedure.  The normal
curve that interests us is the vertical curve, and the log-parameter
along it has been changed by $2\pi i$.

Since
the log-parameter along the normal curve is a multiple of $2\pi
i$ 
(Corollary \ref{cor:logtoparity}), we can apply a multiple of the
above procedure to make sure that desired log-parameter along the
curve is zero before we perform the blow-down.

\begin{figure}[ht!]
  \centering
    \includegraphics[width=.8\hsize]{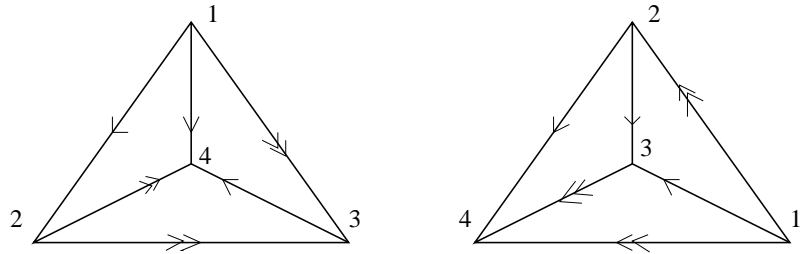}
 \caption{A quasi-simplicial ordered $3$--cycle}
\label{fig:fig8}
\end{figure}

This would complete the proof of Theorem \ref{beta exists} except that
we have carried out the proof only for ordered $3$--cycles whose
vertex-orderings come from a global ordering of the $0$--simplices of
$K$.  We now show that this stronger condition is unnecessary. We will
also show, as promised earlier, that $K$ need only be a
quasi-simplicial complex (closed simplices do not necessarily embed in
$|K|$). An example of such a quasi-simplicial ordered $3$--cycle is
given by gluing the two $3$--simplices in Figure~\ref{fig:fig8} by
matching each face of the left $3$--simplex with the correspondingly
decorated face of the right $3$--simplex (the decorations in question
are the single and double arrows on the edges that bound the face).
The arrows on the edges order the vertices of each $3$--simplex as
shown, and, by construction, these orderings match on common faces.
The resulting complex has one $0$--simplex, two $1$--simplices, four
$2$--simplices, and two $3$--simplices. (This $3$--cycle is simply
connected, so it cannot give a non-trivial element of $H_3(G;\Z)$ in
Theorem \ref{beta exists}, but we will return to it when discussing
invariants of $3$--manifolds, since $|K|-K^{(0)}$ is Thurston's ideal
triangulation of the figure eight knot complement.)
\begin{proposition}
  \label{prop:local order} Any quasi-simplicial ordered $3$--cycle $K$
  is related by Pachner moves that respect vertex-orderings of
  $3$--simplices to a simplicial $3$--cycle $K'$ whose vertex-orderings
  are induced by a global ordering of the $0$--simplices of $K'$.
\end{proposition}
\begin{proof}
  We will start with an arbitrary ordering of the $0$--simplices of
  $K$. At first this ordering will not induce the given
  vertex-orderings of simplices. We will perform various Pachner moves
  on $K$. Each time we perform a $1\to 4$ move we will add the new
  $0$--simplex as a new maximal element both in our ordering of
  $0$--simplices and also in the vertex-orderings of the new
  $3$--simplices created by the move.  At each stage in the process,
  the ordering of the two vertices of a $1$--simplex will be induced by
  the global ordering of the $0$--simplices except maybe if both
  vertices are ``old'' $0$--simplices (ie, they occurred in $K$).  We
  will therefore be done once no two old $0$--simplices are joined by a
  $1$--simplex.
  
  We will assume no $3$--simplex of $K$ meets itself across a $2$--face
  (do a $1\to 4$ Pachner move on any such simplex if it occurs).
  Suppose $K$ has $n$ $3$--simplices. We first do $1\to 4$ Pachner
  moves on each of the $3$--simplices of $K$, creating $n$ new
  $0$--simplices. Next, for each old $2$--simplex (ie, a $2$--simplex
  that was already in $K$) we do a $2\to 3$ Pachner move on the two
  $3$--simplices that meet across it. Each $3$--simplex of the resulting
  complex is in the star of a unique ``old'' $1$--simplex. So if we
  operate in the star of an old $1$--simplex it will not affect the
  star of any other old simplex. 
  
  If there are three or more three $3$--simplices in the star of an old
  $1$--simplex, we reduce the number to three by repeatedly doing
  $2\to3$ moves on pairs of adjacent simplices in its star, and then
  finally do a $3\to 2$ move on the three $3$--simplices in its star to
  remove the old $1$--simplex.  If an old $1$--simplex has just two
  $3$--simplices in its star we do a $1\to4$ move on one of them and
  then do a $3\to 2$ move to remove the old $1$--simplex.
  
  In this way we alter the complex by Pachner moves until all the old
  $1$--simplices have been removed. Then no two old $0$--simplices are
  still joined by a $1$--simplex, so the proof is complete. Since
  Pachner moves do not affect the represented element in
  $\ebloch(\C)$, this also completes the proof of Theorem \ref{beta
    exists}.
\end{proof}

\medskip
The parity condition in the definition of flattening is probably
essential to Theorem \ref{beta exists}, but it's failure can at most
change the resulting element of $\ebloch(\C)$ by $2$--torsion:

\begin{lemma}\label{le:noparity}
  If we represent an element $\alpha\in H_3(G;\Z)$ as in Theorem
  \ref{beta exists} but without requiring the condition on parity then
  the resulting element differs from $\lambda(\alpha)$ at worst by the
  element of order $2$ in $\ebloch(\C)$.
\end{lemma}
\begin{proof}
  Suppose the log-parameter condition is satisfied but the parity
  condition is not satisfied for the flattening of $K$ and let $\mu$
  be the resulting element in $\ebloch(\C)$. Parity along paths
  determines an element of $H^1(|K|-K^{(0)};\Z/2)$ and hence a
  $2$--fold cover of $K$ (possibly branched at some $0$--simplices). The
  lift of the flattening to this cover will satisfy the parity
  condition and hence represent $\lambda(2\alpha)=2\lambda(\alpha)$.
  But it clearly also represents $2\mu$, so
  $2(\mu-\lambda(\alpha))=0$.
\end{proof}
\section{Cheeger--Chern--Simons class}\label{sec:ccs}
The following is simply a restatement of Theorem
\ref{th:main}.
\begin{theorem}\label{th:ccs}
  The homomorphism $\lambda\co H_3(G;\Z)\to \ebloch(\C)$ is an
  isomorphism and it fits in a commutative square:
$$
\begin{CD}
   H_3(G;\Z) @>\lambda>> \ebloch(\C)\\ @VcVV  @VRVV \\
\C/\pi^2\Z @>>> \C/\pi^2\Z\\
\end{CD}
$$
Here $G$ is, as usual, $\PSL(2,\C)$ with discrete topology, $c$ is
the Cheeger--Chern--Simons class, and $R$ is the Rogers dilogarithm map of
Proposition \ref{R}.
\end{theorem}
\begin{proof} We first prove commutativity of the diagram.
  The imaginary part of both $c$ and $R$ gives volume (this was proved
  for $R$ in slightly different language in \cite{neumann}; see also
  \cite{dupont1}) so the only issue is the real part. As pointed
  out in \cite{dupont1}, to prove that the real part of the above
  diagram commutes, we may replace $\PSL(2,\C)$ by $\PSL(2,\R)$.
  
  For an element of $H_3(\PSL(2,\R))$, if we represent it by an
  ordered $3$--cycle $K$ labeled in $\PSL(2,\R)$ and then choose the
  point $z$ in Theorem \ref{beta exists} in
  $\partial\H^2\subset\partial\H^3$ then all the ideal simplices are
  flat. A flat ideal simplex has a natural flattening in our sense,
  since two of its dihedral angles are already $0$ and the third is
  $\pm \pi$, so we adjust only the latter angle to zero. If we flatten
  $K$ this way, we clearly satisfy the logarithmic flattening
  conditions, but the parity conditions are less clear. By Lemma
  \ref{le:noparity} this affects the element of $\ebloch(\C)$, and
  hence also the value of $R$ on it, at most by $2$--torsion.
  
  With this flattening, the function $R$ for a simplex with real
  cross-ratio parameter $x$ becomes the real Rogers dilogarithm of $x$
  used by Dupont in \cite{dupont1}.  Thus $R\circ\lambda$ is (up to
  the $2$--torsion mentioned above) the cocyle he considers. He shows
  there that $R\circ\lambda$ equals $c$ if we map to
  $\C/\frac{\pi^2}{6}\Z$ rather than $\C/\pi^2\Z$ (since
  $\C/\frac{\pi^2}{6}\Z= (\C/\pi^2\Z)/(\Z/6)$, the possible $2$--torsion
  discrepancy maps to zero). Thus $(R\circ\lambda-c)\co
  H_3(G;\Z)\to\C/\pi^2Z$ has image in
  $\frac{\pi^2}6\Z/\pi^2\Z\cong\Z/6$. On the other hand,
  $H_3(\PSL(2,\C);\Z)$ is known to be a divisible group by
  \cite{dupont-sah}, so it has no non-trivial finite quotient, so
  $(R\circ\lambda-c)=0$. Thus $R\circ\lambda=c$, as was to be proved.
  
  It remains to show that $\lambda$ is an isomorphism.  If we compose
  $\lambda$ with the map to $\bloch(\C)$ then the kernel is the
  torsion subgroup of $H_3(G;\Z)$, so the kernel of $\lambda$ is a
  subgroup of the torsion subgroup $\Q/\Z$ of $H_3(G;Z)$. But it was
  pointed out in \cite{dupont1} and \cite{dupont-sah} that the
  Cheeger--Chern--Simons map $c$ is injective on the torsion subgroup of
  $H_3(G;\Z)$, so it follows that $\lambda$ must also be injective on
  this subgroup. Thus $\lambda$ is injective. Surjectivity follows
  similarly: we need only prove it on the torsion subgroup, and for
  this it suffices to see that $R$ is an isomorphism on the torsion
  subgroups, which was shown in the proof of Proposition
  \ref{prop:5.5}.
\end{proof}
One can give an explicit computation of $R\circ\lambda$ on the torsion
of $H_3(G;\Z)$ (see also \cite{dupont-sah94}). By \cite{dupont-sah},
the torsion subgroup $\Q/\Z\subset H_3(G;\Z)$ is given by the image of
$\mu/\{\pm1\}\to G$, where $\mu$ is the group of roots of unity
mapping to diagonal matrices in $G$. Thus the torsion is generated by
the images of generators of $H_3(\Z/n;\Z)\cong \Z/n$ in $H_3(G;\Z)$.
We can represent a generator of this $\Z/n\subset H_3(G;\Z)$ by the
lens space $L(n,1)$ with the flat $G$--bundle given by mapping a
generator of $\pi_2(L(n,1))=\Z/n$ to the $2\pi/n$--rotation in
$\operatorname{SO}(2)\subset \PSL(2,\R)\subset G$. We can triangulate
$L(n,1)$ by triangulating its standard fundamental domain as in
Figure~\ref{fig:cycle}. If we order the vertices of each simplex in that
figure in the order south, north, west, east, then the resulting
representative for this class $\alpha_n\in H_3(G;\Z)$ is the sum of
homogeneous simplices
$$\sum_0^n\langle h_1,gh_1,g^jh_2,g^{j+1}h_2\rangle,$$
where $g$ is
the $2\pi/n$--rotation and $h_1$ and $h_2$ are any elements of
$\operatorname{SO}(2)$ chosen to fulfill our requirement that the
$G$--labels of each simplex be distinct. Explicit computation then
gives $R(\lambda(\alpha_n))=\pi^2/n$ modulo $\pi^2$ (we omit the
details).  

\section{Unordered simplices}\label{sec:unordered}

As we have shown in Proposition \ref{prop:local order}, our complexes
can be quasi-simplicial (closed simplices not required to embed).  The
triangulations that arise in practice (eg, ideal triangulations of
complete finite volume $3$--manifolds, see section
\ref{sec:3-manifolds}) are often not simplicial. However, a
quasi-simplicial triangulation does not always admit vertex-orderings
of its $3$--simplices that agree on common faces.  Thus, to use such a
triangulation one may either have to subdivide, or do without
vertex-orderings.

In this section we describe how the theory changes if we discard the
orderings of vertices of $3$--simplices in the definition of
$\ebloch(\C)$ in terms of flattened simplices. We show that the
result is to quotient $\ebloch(\C)$ by its cyclic subgroup of order
$6$.

Suppose we have a flattened simplex with parameter $(z;p,q)$ with $z$
in the upper half plane. If we retain its geometry but reorder its
vertices by an even permutation (which preserves orientation) then the
flattening is replaced by one of $$(1/(1-z); q, -1-p-q)\quad\text{or}
\quad(1-1/z;-1-p-q,q),$$ so the element in $\ebloch(\C)$ is
changed by subtracting one of 
$$[z;p,q]-[1/(1-z);q,-1-p-q]\quad\text{or}\quad
[z;p,q]-[1-1/z;-1-p-q,p].$$
Similarly, if we reorder its vertices by
an odd permutation then the element in $\ebloch(\C)$ is changed by
subtracting one of
\begin{gather*}
  [z;p,q]+[1/z; -p, 1+p+q],\quad [z;p,q]+[1-z; -q,
-p],\\ \text{or}\quad[z;p,q]+[z/(z-1); 1+p+q, -q].
\end{gather*}
 
\begin{proposition}
  If $z$ is in the upper half plane then, with $\chi$ as in
  Proposition \ref{prop:5.5} and Theorem \ref{th:commutdiag}
  $$\begin{aligned}{}
    [z;p,q]-[1/(1-z);q,-1-p-q]&=\chi(e^{\pi i/3 +q\pi i})\\
[z;p,q]-[1-1/z;-1-p-q,p]&=\chi(e^{-\pi i/3 +p\pi i})\\
  [z;p,q]+[1/z; -p, 1+p+q]&=\chi(e^{p\pi i})\\
[z;p,q]+[1-z; -q,-p]&=\chi(e^{\pi i/3})\\ 
[z;p,q]+[z/(z-1); 1+p+q, -q]&=\chi(e^{2\pi i/3 +q\pi i})
  \end{aligned}$$
\end{proposition}
\begin{proof}
  The last three of the above equations correspond to permutations
  which exchange two vertices. Since two of these involutions,
  together with the Klein four-group (which leaves flattenings
  unchanged) generate the full symmetric group, it suffices to prove
  the third and fourth equations; the others are then easy
  consequences. We start with the third equation.  Recall that
\begin{equation}\label{eq:5t1}
\begin{matrix} [x;p_0,q_0]-[y;p_1,q_1]+[\frac
  yx;p_1-p_0,q_2]-{}\\{}-[\frac{1-x^{-1}}{1-y^{-1}};p_1-p_0+q_1-q_0,q_2-q_1]+[\frac{1-x}{1-y};q_1-q_0,q_2-q_1-p_0]=0
\end{matrix}\quad
\end{equation}
  is an instance of the five-term relation whenever all five of $x$,
  $y$, $\frac yx$, $\frac{1-x}{1-y}$, and $\frac{1-x^{-1}}{1-y^{-1}}$
  are in the upper half plane. Similarly one checks that it is an
  instance of the five-term relation if the first three are
  in the lower half plane and the last two in the upper half plane. In
  particular, if all five are in the upper half plane, then replacing
  $x$ by $x^{-1}$ and $y$ by $y^{-1}$ we get a relation
\begin{equation}\label{eq:5t2}
\begin{matrix} [x^{-1};p_0,q_0]-[y^{-1};p_1,q_1]+[\frac
  xy;p_1-p_0,q_2]-{}\\{}-[\frac{1-x}{1-y};p_1-p_0+q_1-q_0,q_2-q_1]
  +[\frac{1-x^{-1}}{1-y^{-1}};q_1-q_0,q_2-q_1-p_0]=0\,.
\end{matrix}\quad
\end{equation}
Putting $p_0=p_1=q_0=q_1=0$ and $q_2=q$ in equation (\ref{eq:5t1}) and
$p_0=p_1=0$, $q_0=q_1=1$, $q_2=q+1$ in (\ref{eq:5t2}) and then adding,
we get:
\begin{equation}
  \label{eq:1}
([x;0,0]+[x^{-1};0,1])-([y;0,0]+[y^{-1};0,1])+([\frac yx;0,q]+
[\frac xy;0,q+1])=0
\end{equation}
This equation holds if $x$, $y$ and $\frac yx$ are suitably positioned
in the upper half plane. It will continue to hold by analytic
continuation if we vary $x$ and $y$ so long as none of $x$, $y$, or
$\frac yx$ strays out of the upper half plane. If we vary $x$ and $y$
in the upper half plane so that $\frac yx$ crosses interval $(0,1)$
and $\frac xy$ therefore crosses $(1,\infty)$ we get the equation
$$([x;0,0]+[x^{-1};0,1])-([y;0,0]+[y^{-1};0,1])+([\frac yx;0,q]+
[\frac xy;0,q-1])=0,
$$
and exchanging $x$ and $y$ in this equation gives
\begin{equation}
  \label{eq:2}
  ([y;0,0]+[y^{-1};0,1])-([x;0,0]+[x^{-1};0,1])+([\frac
yx;0,q-1]+[\frac xy;0,q] )=0,
\end{equation}
Replacing $q$ by $q+1$ in equation (\ref{eq:2}) and then adding to
equation (\ref{eq:1}) gives, with $z:=y/x$,
$$2\bigl([z;0,q]+[z^{-1};0,q+1]\bigr)=0\,.$$
Here $z$ is now arbitrary
in the upper half plane. Thus $[z;0,q]+[z^{-1};0,q+1]$ is in the
cyclic subgroup of order $2$ in $\ebloch(\C)$. On the other hand, we
know by the proof of Proposition \ref{prop:5.5} that $R$ is an
isomorphism on the torsion subgroup of $\ebloch(\C)$ and one checks
easily that $R([z;0,q]+[z^{-1};0,q+1])=0$. Thus
$[z;0,q]+[z^{-1};0,q+1]=0$. As we analytically continue this relation
by letting $z$ go repeatedly around the origin crossing the intervals
$(-\infty,0)$ and $(0,1)$ (so $z^{-1}$ crosses $(-\infty,0)$ and
$(1,\infty)$) we get the relations
$$[z;p,q]+[z^{-1};-p,p+q+1]=0$$
for all even $p$.

Now, using the notation of the proof of
Proposition \ref{prop:5.5} and Lemma \ref{super transfer} we find 
$$\begin{aligned}{}
  [z;-1,q]+[z^{-1};1,q]&=([z;0,q]+[z^{-1};0,q+1])+([z;-1,q]-[z;0,q]) +\\
&\quad +([z^{-1};1,q] -[z^{-1};0,q])+([z^{-1};0,q]-[z^{-1};0,q+1])\\
&=0+\{1-z,0\}-\{1-z^{-1},0\}-\{z,0\}\\
&=\chi(1{-}z)-\chi(1{-}z^{-1})-\chi(z)=\chi\Bigl(\frac{1{-}z}{(1{-}z^{-1})z}\Bigr)
=\chi(-1)\,.
\end{aligned}$$
Analytically continuing this as $z$ circles the origin as before gives
us the equations
$$[z;p,q]+[z^{-1};-p,p+q+1]=\chi(-1)$$
for all odd $p$. Thus the third equation of the proposition is proved.

For the fourth equation of the proposition, we already know by Lemma
\ref{1-x} that $[z;p,q]+[1-z; -q,-p]=2[1/2;0,0]$, so we need to show
that $2[1/2;0,0]=\chi(e^{\pi i/3})$.  Denote $\omega=e^{\pi i/3}$.
Then $1-\omega=\omega^{-1}$ so
$[\omega;0,0]+[\omega^{-1};0,0]=2[1/2;0,0]$. We have just shown that
$[\omega;0,0]+[\omega^{-1};0,1]=0$. And by definition of $\chi$, we
have
$[\omega^{-1};0,0]-[\omega^{-1};0,1]=-\chi(\omega^{-1})=\chi(\omega)$.
Adding the last two equations gives
$[\omega;0,0]+[\omega^{-1};0,0]=\chi(\omega)$, whence
$2[1/2;0,0]=\chi(\omega)$, as desired (in fact, $[1/2;0,0]=\chi(e^{\pi
  i/6})$, as can be seen now by applying $R$).
\end{proof}
\begin{corollary}
  If we define groups analogous to $\eprebloch(\C)$ and $\ebloch(\C)$
  using flattened simplices modulo the lifted five term relation, but
  ignoring orderings of vertices, we obtain the quotients
  $\eprebloch(\C)/C_6$ and $\ebloch(\C)/C_6$, where $C_6$ is the
  unique cyclic subgroup of order $6$.
\end{corollary}

\section{Invariants of $3$--manifolds}\label{sec:3-manifolds}

Suppose we have a compact oriented hyperbolic $3$--manifold $M$, or,
more generally, a compact oriented manifold $3$--manifold with a flat
$G$--bundle (equivalently, a homomorphism $\pi_1(M)\to G$; as usual $G$
denotes $\PSL(2,\C)$). Then Theorems \ref{beta exists} and
\ref{th:main} give a computation of the corresponding homology class
in $H_3(G;\Z)$ and its Cheeger--Chern--Simons class from a triangulation of
$M$.

In this section we extend this to the case that $M$ is an oriented
hyperbolic manifold which is complete of finite volume but not
necessarily compact. We will also extend to the sort of ideal
triangulations that the programs Snap \cite{snap, coulson-et-al} and
Snappea \cite{snappea} use for compact hyperbolic manifolds. The
underlying complex of such a ``Dehn filling triangulation'' is
homeomorphic not to $M$ itself, but to the result of collapsing to a
point a simple closed curve of $M$, so Theorem \ref{beta exists} does
not apply directly.

Suppose $M$ is a non-compact oriented complete hyperbolic manifold of 
finite volume. It is known by Epstein
and Penner \cite{epstein-penner} that $M$ has a triangulation by ideal
hyperbolic polytopes, and by subdividing these polytopes we may obtain
a triangulation by ideal simplices. As has often been pointed out
(eg, \cite{neumann-yang}), after subdividing the polytopes into
ideal simplices, the subdivisions may not match across faces of the
polytopes. One can mediate between non-matching subdivisions by
including flat ideal simplices in the triangulation.  It is still
unknown if this is actually necessary in any example.  It has been
shown by Petronio and Weeks \cite{petronio-weeks} that flat ideal
simplices are usually not a serious issue.  This is so for our
arguments, in fact for us even ``folded back tetrahedra'' are allowed.

In \cite{neumann} it is shown that any ideal triangulation of a
hyperbolic $3$--manifold has a flattening in the sense of Theorem
\ref{beta exists}. If we can appropriately order vertices then the
arguments of this paper imply that we get an element of
$\ebloch(\C)$.  In general we cannot do this, so, by the previous
section, we only get an element of $\ebloch(\C)/C_6$.
We can refine our triangulation to make sure that we can order the
$0$--simplices. The triangulations we must consider for this are hybrids
of usual triangulations and ideal ones.
 
Suppose $M$ has $h$ cusps, and consider the end compactification of
$M$, obtained by adding points $p_i$, $i=1,\dots,h$, to compactify
each cusp. If $K$ is a complex that triangulates this compactification
of $M$ then $p_1,\dots,p_h$ will be $0$--simplices of $K$. We consider
$K_0:=K-\{p_1,\dots,p_h\}$ to be a ``triangulation'' of the
non-compact manifold $M$. Some of the simplices of $K$ have one or
more vertices among $\{p_1,\dots,p_h\}$; we call these ``ideal
vertices.''

We wish to assign ideal simplex shapes to all the simplices of $K$ in
an appropriate fashion. This is maybe most easily visualized as
follows. We position our triangulation of $M=\H^3/\Gamma$ so that it
triangulates $M$ by geodesic hyperbolic simplices, with vertices
$p_i$ ideal. We then lift to a triangulation of $\tilde M=\H^3$.
Finally, we move each non-ideal vertex to $\partial\overline \H^3$ in
a $\Gamma$--equivariant fashion, making sure that the resulting ideal
simplices are non-degenerate (all four vertices distinct).

\begin{remark}
  We can also describe this is terms of a $G$--valued $1$--cocycle
  ``{relative to the cusps}.''  We choose a lift $\tilde p_i\in \H^3$
  for each $p_i$ and denote by $P_i$ the parabolic subgroup of $G$
  that fixes $\tilde p_i$.  The $G$--valued $1$--cocycle relative to the
  cusps is a $1$--cocycle in the usual sense except that its value on
  an edge that ends (resp.\ starts) in an ideal vertex $p_i$ is only
  well defined modulo right (resp.\ left) multiplication by elements
  of $P_i$. This $1$--cocycle is well defined up to the usual action of
  $0$--cocycles, but only $0$--cocycles that take value $1\in G$ on
  ideal vertices are permitted. The $1$--cocycle depends on the choice
  of the lifts $\tilde p_i$; a change of lift changes the cocycle on
  edges to $p_i$ by right multiplication by the element of $G$ that
  moves the new lift to the old.

There is then an equivalent $G$--labeling of the $3$--simplices as in
Section \ref{subsec:cycles}, except that the label of an ideal vertex
is only defined up to right multiplication by elements of the
corresponding parabolic subgroup.  We then assign ideal simplex shapes
to the simplices of $K$ by choosing a point
$z\in\partial\overline\H^3$ for each $0$--simplex of $K$, as in Remark
\ref{rem:choosing z}, but with the restriction that the point $z$
chosen for an ideal vertex is the fixed point $\tilde p_i$ of the
corresponding parabolic subgroup. A simplex with $G$--labels
$g_1,g_2,g_3,g_4$ and whose vertices have been assigned points
$z_1,\dots, z_4$ in $\partial\overline\H^3$ receives the ideal simplex
shape $\langle g_1z_1,\dots,g_4z_4\rangle$.
\end{remark}

Suppose we have triangulated $M$ and assigned ideal simplex shapes to
the simplices as above.  We also assume that we have chosen the
resulting complex $K$ so
that it admits orderings of its $3$--simplices that agree on common
faces, and we fix such an ordering.  The arguments of sections
\ref{sec:proof} and \ref{sec:proof2} then go through to show the
following theorem:
\begin{theorem}\label{th:general hyperbolic}
  There exists a strong flattening of $K$ as in Theorem
  \ref{flattening exists}.  The resulting element
  $\sum_i\epsilon_{i}[x_i;p_i,q_i]\in\eprebloch(\C)$ is in
  $\ebloch(\C)$ and only depends on the hyperbolic manifold $M$. We
  denote it $\hat\beta(M)$.\qed
\end{theorem}

Since $\ebloch(\C)\cong H_3(G;\Z)$, this gives a
``fundamental class'' in $H_3(G;\Z)$ for any oriented
complete finite volume hyperbolic $3$--manifold (as usual,
$G=\PSL(2,\C)$ with discrete topology). 

We pointed out in the Introduction that this class is easy to define
algebraically if $M$ is compact. I am grateful to the referee for a
comment that led me to a similar elementary description of this class
if $M$ has cusps.  $M$ can be compactified to a manifold $\bar M$ with
$\partial \bar M$ consisting of tori. The fundamental class $[\bar
M,\partial \bar M]\in H_3(\bar M,\partial\bar M)$ determines a class
in $\beta(\bar M,\partial\bar M)\in H_3(BG, BP;\Z)$, where
$P=\left\{\left({1~*\atop0~1}\right)\right\}\subset\PSL(2,\C)$ (as
usual, all groups are considered with their discrete topology). The
long exact sequence for the pair $(BG,BP)$ simplifies to
$$0\to H_i(BG;\Z)\to H_i(BG,BP;\Z)\to H_{i-1}(BP;\Z)\to
0\,,$$
for $i>1$, 
since the
inclusion $P\to G$ factors through the Borel subgroup $B
=\left\{\left({*~*\atop0~*}\right)\right\}\subset\PSL(2,\C)$ and
$H_i(BP;\Z)\to H_i(BB;\Z)$ is zero for $i>0$. 
\begin{proposition}
  There is a
natural splitting
$$\rho\co H_i(BG,BP;\Z)\to H_i(BG;\Z)$$ of the
above sequence for $i>1$.
\end{proposition}
\begin{proof}
  Denote 
\begin{equation*}
  P_2:=\left\{
\begin{pmatrix}
  2^k&b\\0&2^{-k}
\end{pmatrix}
~|~ k\in\Z, b\in \C\right\}\,.
\end{equation*}
One computes that $H_i(BP_2;\Z)=0$ for $i>0$. Thus the diagram
\begin{equation*}
  \begin{CD}
0 @>>> H_i(BG;\Z) @>>> H_i(BG,BP;\Z) @>>> H_{i-1}(BP;\Z) @>>> 0 \\
&& &&  @VVV @|  \\
&& 0 @>>> H_i(BP_2,BP;\Z)  @>\cong >> H_{i-1}(BP;\Z) @>>> 0
\end{CD}
\end{equation*}
splits the top sequence. Note that the splitting $H_{i-1}(BP)\to
H_i(BG,BP)$ is not canonical, since it depends on the ``2'' in the
definition of $P_2$, but it is not hard to see that its image does not
depend on this choice, so the resulting splitting
$$\rho\co H_i(BG,BP;\Z)\to H_i(BG;\Z)$$
is canonical.
\end{proof}
We omit the proof that $\rho(\beta(\bar M,\partial\bar M))$
is indeed the element
$\hat\beta(M)$ of Theorem \ref{th:general hyperbolic}.


\begin{example}\label{ex:fig8}
  Figure~\ref{fig:fig8} gives a well-known ideal triangulation for the
  complete hyperbolic structure on the figure eight knot complement
  (see, eg, \cite{thurston, neumann-zagier}). This manifold is
  denoted $m004$ in the Callahan-Hildebrand-Weeks cusped census (see
  \cite{census,snappea,snap}). Both simplices are regular ideal
  simplices but the vertex-order of the second simplex is opposite to
  its orientation. The cross-ratio parameters are therefore
  $\omega:=e^{\pi i/3}$ and $\omega^{-1}$ respectively, and
  $\hat\beta(m004)=[\omega;p,q]-[\omega^{-1};r,s]$ for suitable
  $p,q,r,s$. In the next section we will see that we get a flattening
  of the complex $K$ if and only if $q=-1-2p$ and $(r,s)=-(p,q)$.
  Choosing, eg, $p=0$, we see
  $$\hat\beta(m004)=
  [\omega;0,-1]-[\omega^{-1};0,1].$$
  If $$V_0=1.014941606409653625021202554\dots$$ denotes the volume of the
  regular ideal tetrahedron, then
  $$
  R(\omega;p,q)=\frac{(2p-2q-1)\pi^2}{12} + V_0i \,,$$
  $$
  R(\omega^{-1};-p,-q)=\frac{(2p-2q-1)\pi^2}{12}- V_0 i\,.$$
  Thus
$$(\vol+i\cs)(m004)=
2V_0\,,$$
recovering the known volume and Chern--Simons invariant of the
figure eight knot complement.
\end{example}

We now describe what happens when we deform our hyperbolic structure
by hyperbolic Dehn surgery to perform a Dehn filling on $M$ (see,
eg, \cite{thurston, neumann-zagier}). Topologically, the Dehn filled
manifold $M'$ differs from $M$ in that some of the cusps have been
filled by solid tori, which adds new closed geodesics at these cusps
(the cores of the solid tori). Let $\lambda_j$ be the complex length
of the geodesic (length plus $i$ times torsion) added at the $j$-th
cusp. If no geodesic has been added at the $j$-th cusp we put
$\lambda_j=0$.

The deformation deforms the representation $\Gamma=\pi_1(M)\to G$.  If
we keep the points $\tilde p_i$ at fixed points of the images of the
cusp subgroups of $\Gamma$ then our triangulation deforms to what is
called a ``degree one ideal triangulation'' in \cite{neumann-yang} (it
is not a genuine triangulation because the topology of $K$ has not
changed, so its topology is not the topology of $M'$). By deforming
the flattened ideal simplex parameters continuously we obtain
flattened ideal simplex parameters $[x_i';p'_i,q'_i]$ after
deformation. 
\begin{theorem}\label{th:conj true}
  With the above notation, the element in $\ebloch(\C)$ represented by
  the Dehn filled manifold $M'$ is
  $$\hat\beta(M')=-\sum_{j=1}^h
  \chi(e^{\lambda_j})+\sum_i\epsilon_i[x'_i;p'_i,q'_i]$$
with $\chi$ as in Theorem \ref{th:commutdiag}. Moreover, 
$$(\vol +i\cs)(M')=-\frac{\pi }2\sum_{j=1}^h \lambda_j
-i\sum_i\epsilon_iR(x'_i;p'_i,q'_i)$$
\end{theorem}
The second formula of this theorem was proved up to a constant
(depending on $K$ but independent of the Dehn filling) in
\cite{neumann}.  The constant was conjectured to be a multiple of
$i\pi^2/6$. Since \cite{neumann} used unordered simplices, the
constant in that version is indeed a multiple of $i\pi^2/6$.  Versions
prior to 1.10.2 of the program Snap \cite{snap, coulson-et-al} use this
formula but ignores the parity condition, giving answers accurate to
$i\pi^2/12$. Snappea \cite{snappea} uses a version, also from
\cite{neumann}, that is accurate to an unknown constant. Both programs
then bootstrap this to an accurate computation of the Chern--Simons
invariant for any manifold that can be related by a sequence of Dehn
drillings and fillings to one with known Chern--Simons invariant. The
value they print is $\cs(M)/2\pi^2$, hence well defined modulo $1/2$.

Snap now uses the formulae of this paper and can compute Chern--Simons
for any manifold. For compact manifolds accessible by the
bootstrapping method, Snap also computes the eta-invariant using a
related formula of \cite{neumann-meyerhoff}; from this the Riemannian
Chern--Simons (defined modulo $2\pi^2$, hence modulo $1$ in their
normalization) is then computed. 

\begin{proof}[Proof of Theorem \ref{th:conj true}]
  For simplicity of
exposition we assume $M$ has just one cusp represented by the ideal
vertex $p$ of $K$.

The proof depends on the fact that the deformed flattening  parameters
$[x_i';p_i',q_i']$ still satisfy our flattening conditions
around edges, but there is a natural flattening
condition at the filled cusps that is not satisfied, which leads to the
``correction term'' $\sum_{j=1}^h
  \chi(e^{\lambda_j})$ in the theorem. 

The link $L_p$ of vertex $p$ of $K$ is a torus. Let $C$ be a
simplicial curve in $L_p$ that is a meridian of the solid torus that
the cusp has been replaced by.  We can obtain a genuine triangulation
of $M'$ by the blow-up procedure of Section \ref{sec:proof3}, blowing
up using the curve $C$. As in that section, we want to extend the
flattening after blowing up so that matching pairs of added simplices
have the same flattening and therefore cancel in the computation of
$\hat\beta(M')$. To be able to extend the flattening of $K$ after
blowing up, we need that the log-parameter along a normal curve
parallel to $C$ is zero. As it stands, it is not zero, since it
represents the logarithm of the derivative of the holonomy of the
meridian curve.  It was zero before deformation, but after deformation
the meridian curve represents a full rotation about the core curve of
the solid torus added by the Dehn filling, so the logarithm of the
derivative of its holonomy is $2\pi i$ (see, eg, \cite{thurston,
  neumann-zagier} for more detail).  We can use the procedure of
Section \ref{sec:proof3} to correct this: we must modify flattenings
on the simplices traversed by a normal curve representing a longitude
of the added solid torus as indicated in Figure~\ref{fig:annulus} (the
longitude runs horizontally across the center of the figure). The
effect on each modified simplex is to change its contribution to
$\sum_i\epsilon_i[x_i;p_i,q_i]$ by something of the form
$[x;p,q\pm1]-[x;p,q]$ or the equivalent after a permutation of
vertices. By the calculations of Proposition \ref{prop:5.5} this is
$\pm\chi(z)$, where $z$ is the cross-ratio parameter corresponding to
the edge of the $3$--simplex that the curve passes (so $z$ is one of
$x$, $1/(1-x)$, $1/x$, etc.). As we sum these contributions over the
affected simplices with appropriate signs, they sum to $-\chi(L)$,
where $L$ is product of the corresponding cross-ratio parameters or
their inverses for the edges passed by the normal curve. This $L$ is
the derivative of holonomy of the longitude and $\log(L)$ is the
complex length $\lambda$ of the core curve of the solid torus added by
the Dehn filling (see, eg, \cite{thurston, neumann-zagier}). This
proves the first equation of the theorem.

For the second equation we use the fact that
$R(\hat\beta(M'))=i(\vol+i\cs)(M')$. We know this if $M'$ is
compact and we will deduce it below if  $M'$ is non-compact, so we
assume it for now. $R\circ\chi\co \C^*\to \C/\pi^2\Z$ is the
map $z\mapsto \frac{\pi i}2\log z$, so $R\chi(e^{\lambda})=\frac{\pi
  i}2\lambda$. Applying $\frac1iR$ to the first equation of the
theorem thus gives the second.
\end{proof}

Meyerhoff's extension of Chern--Simons invariant to
cusped hyperbolic manifolds in \cite{meyerhoff} is given by defining
$(\vol+i\cs)(M)$ to be the limit of $(\vol+i\cs)(M')$\break $+\frac\pi
2\sum_{j=1}^h\lambda_j $ as $M'$ approaches $M$ in hyperbolic Dehn
surgery space. The case that $M'$ is compact in the above Theorem thus
implies:
\begin{corollary}
  For any oriented complete finite volume hyperbolic $3$--mani\-fold $M$
  $$i(\vol+i\cs)(M)=R(\hat\beta(M))\in\C/\pi^2\Z,$$
where $cs(M)$ is Meyerhoff's extension of the Chern--Simons invariant if
$M$ is non-compact.\qed
\end{corollary}

The ``correction term'' $-\sum_{j=1}^h \chi(e^{\lambda_j})$ in Theorem
\ref{th:conj true} arises because the sum of flattening parameters
corresponding to a meridian of the a solid torus added at a cusp by
Dehn filling is $2\pi i$ rather than zero. The proof of the theorem
shows that flattening parameters can be chosen so that this sum is
zero, and then the correction term is not there. This gives an analog
of Theorem \ref{beta exists} for these triangulations:
\begin{theorem}\label{th:deg one}
  Consider the degree one triangulation of $M'$ of
  Theorem \ref{th:conj true}. Then there exists flattenings
  $[x_i';p_i'',q_i'']$ of the ideal simplices which satisfy the
  conditions
\begin{itemize}
\item parity along normal paths is zero; 
\item log-parameter about each edge is zero;
\item log-parameter along any normal path in the neighborhood of a
  $0$--simplex that represents an unfilled cusp is zero;
\item log-parameter along a normal path in the neighborhood of a
  $0$--simplex that represents a filled cusp is zero if the path 
 is null-homotopic in the added solid torus.
\end{itemize}
For any such choice of flattenings we have
$$\hat\beta(M')=\sum_i[x'_i;p''_i,q''_i]\in\ebloch(\C)$$
so $$
i(\vol+i\cs)(M')=\sum_iR(x'_i;p''_i,q''_i)\,.\eqno{\qed}$$
\end{theorem}

\section{Example}

We return to the figure 8 knot complement $m004$ of Example
\ref{ex:fig8} to illustrate the above formulae.  Denote the
cross-ratio parameters of the two simplices in Figure~\ref{fig:fig8} by
$z$ and $y$, and choose flattenings $[z;p,q]$ and $[y,r,s]$. The
consistency condition about
each edge can be read off from the figure. The
condition for the edge labeled with a single arrow is:
\begin{equation}
  \label{eq:5}
  2\log z+\log z'-\log y'-2\log y''=2\pi i
\end{equation}
and the condition for
the other edge turns out to be equivalent to this.  The log-parameters
for the flattenings are
$$\begin{aligned}
  &(\log z+p\pi i,~\log
z'+q\pi i, ~\log z''+(-1-p-q)\pi i),\\
&(\log y+r\pi i,~\log y'+s\pi
i, ~\log y''+(1-r-s)\pi i)\,.
\end{aligned}$$
\begin{figure}[ht!]
  \centering
  \includegraphics[width=.95\hsize,height=.475\hsize]{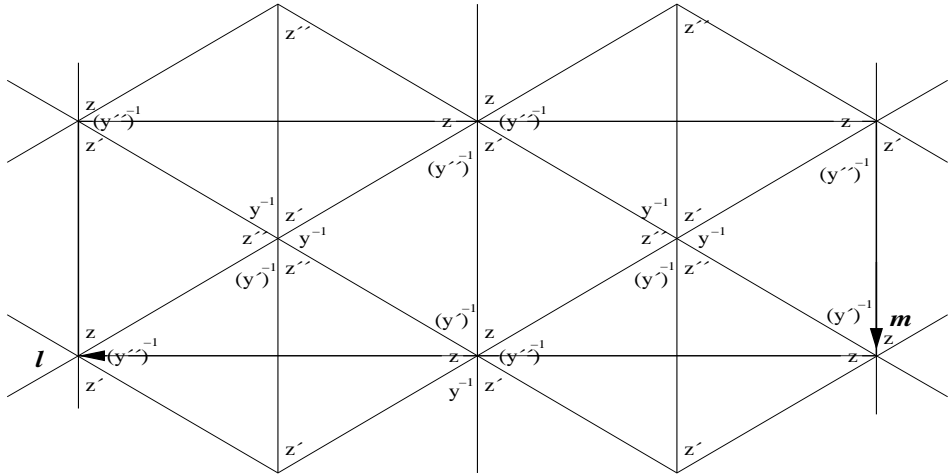}
  \caption{Fundamental domain of the cusp torus, triangulated
  as the link of the $0$--simplex of the ideal triangulation of $m004$.}
  \label{fig:cusp}
\end{figure}
Thus the flattening condition about the edge is
\begin{equation}
  \label{eq:3}
  2p+q-s-2(1-r-s)=-2\,. 
\end{equation}
The log holonomy for a meridian and longitude can be read off from the
fundamental domain of the cusp torus illustrated in
Figure~\ref{fig:cusp} (taken from \cite{neumann-zagier} but modified
since orientation conventions differ; in particular, $w,w',w''$ there
are our $(y'')^{-1},(y')^{-1},y^{-1}$ because of different
vertex-ordering for the second simplex). They are:
\begin{equation}
  \label{eq:7}
u:=\log z''+\log y''\qquad  v:=2\log z-2\log z''\,.
\end{equation}
Their vanishing thus gives flattening conditions
$$\begin{aligned}
(1-r-s)+(-1-p-q)&=0\\
  2p-2(-1-p-q)&=0\,
\end{aligned}$$
which simplify to
\begin{align}
  \label{eq:4}
-p-q-r-s&=0\\
\label{eq:4bis}
  2+4p+2q&=0\,.
\end{align}
Solving equations (\ref{eq:3}), (\ref{eq:4}), (\ref{eq:4bis}) gives 
\begin{equation}
  \label{eq:6}
q=-1-2p,\quad (r,s)=-(p,q).
\end{equation}
Simultaneously solving the consistency condition (\ref{eq:5})
and the cusp conditions $u=v=0$ (see (\ref{eq:7})) gives $z=\omega$,
$y=\omega^{-1}$.  Thus, as promised earlier, we see that 
$\hat\beta(m005)=[\omega;p,q]-[\omega^{-1};-p,-q]$ with $q=-1-2p$. 
Choosing $p=0$, we get:
$$\hat\beta(m005)=[\omega;0,-1]-[\omega^{-1};0,1].$$
Now suppose we deform our parameters to perform hyperbolic
$(\alpha,\beta)$ Dehn filling on this manifold. Let $z,y$ now denote
the deformed values of the simplex cross-ratios (they simultaneously
satisfy the consistency condition (\ref{eq:5}) and the Dehn filling
condition $\alpha u+\beta v=2\pi i$).  Choose $\gamma, \delta\in\Z$ so
$$\left|
  \begin{matrix}
    \alpha&\beta\\ \gamma&\delta
  \end{matrix}\right|=1.$$
As described in \cite{neumann-zagier}, the complex length of the added
geodesic is $\lambda:=-(\gamma u+\delta v)$, so Theorem  \ref{th:conj
  true} gives the formulae
\begin{equation}\label{eq:8a}\begin{aligned}
\hat\beta(m004(\alpha,\beta))&=
-\chi(e^{-\gamma u-\delta v})+[z;0,-1]-[y;0,1]\\
(\vol+i\cs)(m004(\alpha,\beta))&=\frac\pi2(\gamma u+\delta
v)-i(R(z;0,-1)-R(y;0,1))\,,
\end{aligned}
\end{equation}
with $u,v$ given by (\ref{eq:7}).

On the other hand, if we use Theorem \ref{th:deg one}, 
the Dehn filling condition $\alpha u+\beta v=2\pi i$ gives the
flattening condition
$$\alpha(-p-q-r-s)+\beta(2+4p+2q)=-2$$
while the parity condition says
$p+q+r+s$ must be even.  A simultaneous solution of these and equation
(\ref{eq:3}) is:
$$q=\gamma-1-2p,\quad r=-2\delta-p,\quad s=-\gamma+4\delta+1+2p\,.$$
Choosing, eg, $p=0$, we get
\begin{equation}
  \label{eq:8}
  \begin{aligned}
\hat\beta(m004(\alpha,\beta))&=[z;0,\gamma-1]-[y;-2\delta,
-\gamma+4\delta+1]\,,\\
(\vol+i\cs)(m004(\alpha,\beta))&=
-i(R(z;0,\gamma-1)-R(y;-2\delta,-\gamma+4\delta+1))\,.
  \end{aligned}
\end{equation}
It is not hard to verify directly the equivalence of 
(\ref{eq:8a}) and (\ref{eq:8}).

\section{Other fields}
Although we have worked over $\C$
in this paper, most of what we have done works for a subfield
$K\subset \C$ if we replace $G=\PSL(2,\C)=\PGL(2,\C)$ by
$G=\PGL(2,K)$. In Theorem \ref{th:commutdiag} we must replace the
group $\mu$ of roots of unity by the group $\mu_K$ of roots of unity
in $K$.  The homomorphism $\lambda\co H_3(\PGL(2,K))\to
\ebloch(K)$ is still defined.
For a hyperbolic $3$--manifold $\H^3/\Gamma$ with
$\Gamma\subset \PGL(2,K)$ the element $\hat\beta(M)$ naturally lies in
$\ebloch(K)$, but we have numerical evidence that it lies in
$\ebloch(k)$, where $k$ is the invariant trace field of $M$ ($k$ is
always contained in $K$ and is generally smaller). The arguments of
\cite{neumann-yang} show that $\hat\beta(M)\in\ebloch(k)$ if $M$ is
non-compact, while some power of 2 times $\hat\beta(M)$ lies in
$\ebloch(k)$ in general.

\end{document}